\documentclass[10pt,a4paper]{article}
\usepackage{t1enc}
\usepackage[utf8]{inputenc}
\usepackage[english]{babel}
\usepackage{graphicx}
\usepackage{amsmath}
\usepackage{amssymb}
\usepackage{mathrsfs}
\usepackage{theorem}
\usepackage[all]{xy}
\usepackage{eufrak}
\usepackage{enumerate}

\newtheorem{thm}{Theorem}[section]
\newtheorem{prop}[thm]{Proposition}

\newtheorem{dfn}[thm]{Definition}
\newtheorem{cor}[thm]{Corollary}
\newtheorem{lem}[thm]{Lemma}

\DeclareMathOperator{\Aut}{Aut}

\DeclareMathOperator{\Gal}{Gal}

\DeclareMathOperator{\Spec}{Spec}
\DeclareMathOperator{\Spf}{Spf}
\DeclareMathOperator{\temp}{temp}
\DeclareMathOperator{\gf}{\pi_1}

\DeclareMathOperator{\ga}{\pi_1^{alg}}
\DeclareMathOperator{\gtop}{\pi_1^{top}}

\DeclareMathOperator{\gt}{\pi_1^{t}}

\DeclareMathOperator{\gorb}{\pi_1^{orb}}
\DeclareMathOperator{\gtemp}{\pi_1^{temp}}

\DeclareMathOperator{\Covalg}{Cov^{alg}}
\DeclareMathOperator{\Covtop}{Cov^{top}}
\DeclareMathOperator{\Covtemp}{Cov^{temp}}

\DeclareMathOperator{\Ba}{\mcal B^{alg}}

\DeclareMathOperator{\Bcov}{\mathcal B^{cov}}
\DeclareMathOperator{\Btemp}{\mathcal B^{temp}}

\DeclareMathOperator{\Ker}{Ker}

\DeclareMathOperator{\ab}{ab}
\DeclareMathOperator{\et}{\text{ét}}

\newcommand{\Ens}{\text{Ens}}
\DeclareMathOperator{\Set}{Set}
\DeclareMathOperator{\Stab}{Stab}

\DeclareMathOperator{\Ind}{Ind}

\DeclareMathOperator{\an}{an}
\DeclareMathOperator{\Image}{Im}

\DeclareMathOperator{\sq}{\square}

\DeclareMathOperator{\op}{op}

\DeclareMathOperator{\diametre}{diam}

\newcommand{\findem}{\vspace{-0.25in} \begin{flushright}$\square$\end{flushright}}
\newcommand{\dem}{\emph{\underline{Proof} :\;}}
\newcommand{\da}{\begin{displaystyle}}
\newcommand{\db}{\end{displaystyle}}

\newcommand{\mcal}{\mathcal}
\newcommand{\mbf}{\mathbf}

\newcommand{\mbb}{\mathbb}
\newcommand{\fk}{\mathfrak}

\newcommand\sommet{*-[o]{\bullet}}
\begin{document}
\title{Tempered fundamental group and metric graph of a Mumford curve}
\author{Emmanuel Lepage}
\maketitle
\tableofcontents

\section*{Introduction}
This paper is an attempt to give some general results on the
tempered fundamental group of $p$-adic smooth algebraic varieties.\\
The tempered fundamental group was introduced in~\cite[part III]{andre1} as a sort of analog of the topological fundamental group of complex algebraic varieties; its profinite completion coincides with Grothendieck's algebraic fundamental group, but it has itself many infinite discrete quotients in general.\\
Since the analytification (in the sense of Berkovich or of rigid geometry) of a finite étale covering of a $p$-adic varieties is not necessarily a topological covering, André had to use a slightly larger notion of covering.
He defines tempered coverings, which are \'etale coverings in
the sense of de Jong (that is to say that locally on the Berkovich
topology, it is a direct sum of finite coverings) such that, after pulling
back by some finite \'etale covering, they become topological coverings (for the
Berkovich topology). Then the tempered fundamental group is a prodiscrete group that classifies those tempered coverings.  To give a more handful description, if
one has a sequence of pointed finite coverings $((S_i,s_i))_{i\in \mbf N}$ such that the corresponding pointed pro-covering of $(X,x)$ is the
universal pro-covering of $(X,x)$, and if $(S^{\infty}_i,s^{\infty}_i)$ is a
universal topological covering of $S_i$, the tempered fundamental group of
$X$ can be seen as $\gtemp(X,x)=\varprojlim_i \Gal(S^{\infty}_i/X)$.
Therefore, to understand the tempered fundamental group of a variety, one
mainly has to understand the topological behavior of the finite \'etale
coverings of this variety.\\
There are many differences between the tempered fundamental group in the $p$-adic case and the topological fundamental group in the complex case. First, the tempered fundamental group is not discrete in general. It is also much more difficult to describe explicitly (such a description is available for elliptic curves, or more generally abelian varieties as will be studied in this article). As proved in~\cite{mochi} (and recalled in part~\ref{mochipart}), the tempered fundamental group of a curve depends hugely on the combinatorial structure of its stable reduction (this suggests a geometric anabelian behavior of the tempered fundamental group which as no algebraic or complex counterpart). On the other hand, the tempered fundamental group, like the algebraic fundamental group, is also defined over non algebraically closed nonarchimedean field, and thus it interacts interestingly with Galois theory (which gives it a number theoretical interest). However, in this article, we will only study the geometric tempered fundamental group.
We will prove here some results for the tempered fundamental group which are classical results for the profinite fundamental group or the topological
fundamental group of complex varieties: we will link the abelianized tempered fundamental group of a curve to the tempered fundamental group of its Jacobian variety and we will show a Künneth formula for the product of manifolds.
But, to illustrate that it depends much more on the variety itself than the profinite or the complex fundamental group, we will show that one can recover the metric structure of the graph of the stable model of a Mumford curve from the tempered fundamental group.\\

After recalling the basic results on the tempered fundamental group given in~\cite{andre1} (and deduce from them the birational invariance of the tempered fundamental group in proposition~(\ref{birat})), we will recall  the results of Mochizuki concerning the tempered fundamental group of a curve and the stable reduction of this curve:
whereas the profinite fundamental group or
the topological fundamental group of a
smooth curve of type $(g,n)$ only depends in the complex case on $g$ and $n$, the tempered
fundamental group of a curve depends much more on the curve itself. How much is an unsolved problem, which is the leading thread of the present paper.\\
Mochizuki thus proved in~\cite{mochi} that one could recover from the
tempered fundamental group of a curve the graph of its stable
reduction (theorem~(\ref{mochi})), even in the pro-$(p')$ case (we will denote by
$\gtemp(X)^{(p')}$ the group classifying coverings that became topological
after base change by a finite Galois covering of order prime to $p$).\\
 More precisely, the vertices correspond to the classes of
maximal subgroups of $\gtemp(X)^{(p')}$ (thus called verticial subgroup) and
the edges to non trivial classes of intersection of verticial subgroups.\\

We will then study the tempered fundamental group of an abelian variety in the geometric case. We will be able to give an explicit description of the tempered fundamental group, as was already done
for elliptic curves in~\cite[III.2.3.2]{andre1}. The point is that the finite étale
coverings are well understood (they are abelian varieties, and they are
dominated by the multiplication by $n$ on $A$ for some $n$), and the
topological coverings are well understood too, thanks to the $p$-adic
uniformization of abelian varieties.\\
We will deduce from this that $\gtemp(A)$ is isomorphic to $\widehat{\mbf
  Z}^{2g-d}\times \mbf Z^d$ (see part~\ref{abvar}), where $g$ is the dimension of $A$ and $d$ is the
rank of the topological fundamental group of $A$ (which is also the dimension of the toric part of the semistable reduction).\\
One can then prove that, just as in the profinite case or the topological case,
the abelianized tempered fundamental group of a curve $X$ is canonically isomorphic to
the tempered fundamental group of its Jacobian variety $A$ (theorem~(\ref{jacobian})). The proof combines the birational morphism $X^g/\frak S_g\to A$, the birational invariance of the tempered fundamental group, and the Künneth formula that will be proved in the next part.\\

For more general cases, one does not have such nice descriptions of the
topological behavior of a smooth variety, not to mention the
topological behavior of all its finite etale coverings. Nevertheless,
thanks to de Jong's alteration results (see~\cite{dJ2} and~\cite{dJ3}) and to Berkovich's work
(see~\cite{berk2}, recalled in part~\ref{berkspaces}) on the topological structure of the generic fiber of a
polystable formal scheme over the integral ring of a complete
nonarchimedean field, one can prove that some Zariski open dense subset $U$ of
the smooth scheme $X$ is homotopically equivalent to a certain polysimplicial set given by
the combinatorial structure of the special fiber of some alteration
of the variety.\\
One can then deduce from the results of~\cite{berk2} that for any isommetric embedding of $K$ in an algebraically closed field $K'$, $|U_{K'}|\to |U|$ is a homotopy equivalence. 
Knowing this, and using a proposition of André which says that $U\to X$ induces an isomorphism on topological groups, we will be able to prove the invariance of the
tempered fundamental group of a smooth scheme under algebraically closed extensions (proposition~(\ref{invariance}).\\
In the same way, if $V$ is a Zariski dense open subset of another smooth scheme $Y$ satisfying the same properties as $U$, one can deduce from the results of~\cite{berk2} that $|U\times V|\to |U|\times |V|$ is a homotopy equivalence. Thanks to this, we will prove a Künneth formula for the tempered fundamental group of smooth schemes (proposition~(\ref{prod})).\\ 

Finally, we will prove in theorem~(\ref{mumford}) that, in the case of a Mumford curve, one may in fact recover the length of the edges from
$\gtemp(X)$ (obviously, one needs the whole $\gtemp(X)$ and not only the
pro-$(p')$ version which only depends of the pro-$(p')$ graph of
groups of the stable reduction; thus this can only be done when $p\neq 0$,
contrarily to Mochizuki's result).\\
We will deduce from Mochizuki's results that one can decide from the tempered fundamental group of a curve whether a finite covering is split over some vertex of the curve (and in a similar way, one can decide whether it is ramified over some cusp).\\
Knowing this, we will start by proving this for a punctured projective line and a punctured elliptic curve, to give an insight of the method in simpler cases, where the coverings we will use are described explicitly: for a punctured projective line, we will concentrate on coverings like $(\ )^{p^e}:\mbf G_m\stackrel{z\mapsto z^{p^e}}{\to}\mbf G_m$ (one can in this case easily calculate explicitly the number of preimage of a point of the Berkovich projective line), and for elliptic curves on coverings obtained by patching of those types of coverings.\\
In the more general case of a Mumford curve $X=\Omega/\Gamma$, we will study coverings of the topological uniformization $\Omega$ that descend to some finite topological covering of $X$ that behave like $(\ )^{p^e} $ over big affinoids of $\Omega$. More precisely, every covering of $\Omega$ which descends to $X$ (in some non unique way) is the pullback of $(\ )^{p^e}$ along a $\Gamma$-equivariant invertible section of $\Omega$, and conversely. The equivariance prevents us to simply consider a homography, but we can consider a section of $\Omega$ which is arbitrarily close, on a given affinoid subset of $\Omega$, to an homography. This will be enough to ensure that this covering is split over the same vertices of a big affinoid subset of $\Omega$ as the one obtained by pulling back $(\ )^{p^e}$ along the homography.\\
Using such constructions, one can recover from the tempered fundamental group the length of any loop of any topological covering of $X$. In order to get the whole metric structure of the graph of the stable model of $X$, we will end this article by proving a purely combinatorial result that shows that, if one knows the length of every loop of every covering of a graph whose edges have valency $\geq 3$, one knows the length of every edge (proposition~(\ref{combi})).\\
The recovering of the metric graph of the stable model from the tempered fundamental group can be quite easily extended to the case of punctured Mumford cases, but the proof does not seem to extend easily to more general curves, as one cannot find that easily coverings of order divisible by $p$ where one can say much about the graph of the stable reduction.\\
The proofs here also only consider some very simple coverings (for which we do not even describe completely the graph of the stable model), and it is hard to imagine what can recover from the whole tempered fundamental group.\\

This work is part of a PhD thesis. I would like to thank my advisor, Yves André, for guiding my work, suggesting the main idea of some proofs presented here, and taking the time of reading and correcting this work.

\section{Tempered fundamental group}
\subsection{Definition}
Let $K$ be a complete nonarchimedean field.\\
Following~\cite[§4]{andre2}, a $K$-\emph{manifold} will be a smooth paracompact strictly $K$-analytic space. For example, if $X$ is a smooth algebraic $K$-variety, $X^{\an}$ is a $K$-manifold (and in fact, we will mainly be interested in those spaces). Then, thanks to~\cite{berk2}, such a manifold is locally contractible (we will explain in more detail the results of~\cite{berk2} in section~\ref{berkspaces}).\\
Recall from~\cite{dJ1} that a morphism $f:S'\to S$ is said to be an \emph{étale covering} if $S$ is covered by open subsets $U$ such that $f^{-1}(U)=\coprod V_j$ and $V_j\to U$ is étale finite.\\
Then, André defines tempered coverings as follows:
\begin{dfn} \label{def:rvt:temp}(\emph{\cite[def. 2.1.1]{andre1}})
An étale covering $S' \to S$ is \emph{tempered} if it is a quotient of the composition of a topological covering $T' \to T$ and of a finite étale covering
 $T \to S$.
\end{dfn}
This is equivalent to say that it becomes a topological covering after pullback by some finite étale covering.\\
We denote by $\Covtemp(X)$ the category of tempered coverings of $X$ (with the obvious morphisms).\\
Let $\bar x$ be a geometric point of $X$. Then one has a functor \[F_{\bar x}:\Covtemp(X)\to\Set\] which maps a covering $S\to X$ to the set $S_{\bar x}$. If $\bar x$ and $\bar x'$ are two geometric points, then, according to~\cite[prop 2.9]{dJ1}, $F_{\bar x}$ and $F_{\bar x'}$ are (non canonically) isomorphic.\\
The tempered fundamental group of $X$ pointed at $\bar x$ is \[\gtemp(X,\bar x)=\Aut F_{\bar x}.\]
This is a prodiscrete topological group, for which the stabilizers $(\Stab_{F(S)}(s))_{S\in \Covtemp(X),s\in F_{\bar x}(S)}$ form a basis of open subgroups of $\gtemp(X,\bar x)$. The tempered fundamental group depends on the basepoint only up to inner automorphism (this topological group, considered up to conjugation, will sometimes be denoted simply $\gtemp(X)$).\\
The full subcategory of tempered coverings $S$ for which $F_{\bar x}(S)$ is finite is equivalent to $\Covalg(S)$, so that \[\widehat{\gtemp(X,\bar x)}=\ga(X,\bar x)\] (where $\widehat{\ }$ denotes here, and in the sequel, the profinite completion).\\
For any morphism $X\to Y$, the pullback defines a functor $\Covtemp(Y)\to\Covtemp(X)$. If $\bar x$ is a geometric point of $X$ with image $\bar y$ in $Y$, this gives rise to a continuous homomorphism \[\gtemp(X,\bar x)\to\gtemp(Y,\bar y)\] (hence an outer morphism $\gtemp(X)\to\gtemp(Y)$).\\
One has the analog of the usual Galois correspondence:
\begin{thm}\emph{(\cite[Th. 1.4.5]{andre1})}\label{galcorr} $F_{\bar x}$ induces an equivalence of categories between the category of direct sums of tempered coverings of $X$ and the category $\gtemp(X,\bar x)-\Set$ of discrete sets endowed with a continuous left action of $\gtemp(X,\bar x)$.\end{thm}

If $S$ is a finite Galois covering of $X$, its universal topological covering $S^{\infty}$ is still Galois and every connected tempered covering is dominated by such a Galois tempered covering.\\
If $((S_i,\bar s_i))_{i\in \mbf N}$ is a cofinal projective system (with morphisms $f_{ij}:S_i\to S_j$ which maps $s_i$ to $s_j$ for $i\geq j$) of geometrically pointed Galois finite étale coverings of $(X,\bar x)$, let $((S^{\infty}_i,\bar s^{\infty}_i))_{i\in\mbf N}$ be the projective system, with morphisms $f_{ij}^\infty$ for $i\geq j$, of its pointed universal topological coverings. Then $F_{\bar x}(S^{\infty}_i)=\gtemp(X,\bar x)/\Stab_{F(S^\infty_i)}(\bar s^\infty_i)$ is naturally a quotient group of $\gtemp(X,\bar x)$  for which $s^{\infty}_i$ is the neutral element. Moreover $G$ acts by $G$-automorphisms on $F_{\bar x}(S^{\infty}_i)$ by right translation (and thus on $S^{\infty}_i$ thanks to the Galois correspondence of theorem~(\ref{galcorr})). Thus one gets a morphism $\gtemp(X,\bar x)\to\Gal(S^{\infty}_i/X)$. As $f_{ij}^\infty(s_i^\infty)=s_j^\infty$, these morphisms are compatible with $\Gal(S^{\infty}_i/X)\to \Gal(S^{\infty}_j/X)$.\\
Then, thanks to~\cite[lemma III.2.1.5]{andre1},
\begin{prop} \[\gtemp(X,\bar x)\to\varprojlim \Gal(S^{\infty}_i/X)\] is an isomorphism.\end{prop}

We will also have to use the following result by Andr\'e:
\begin{prop}\label{andre114}\emph{(\cite[prop. III.1.1.4]{andre1})}
Let $\overline S$ be a manifold, and let $Z$ be a Zariski closed nowhere dense reduced analytic subset. Then  any topological covering of $S:=S\backslash Z$ extends uniquely to a topological covering. Thus
$\gtop(S,s)\to\gtop(\overline S,s)$ is an isomorphism.\end{prop}

We can deduce from this that:
\begin{prop}\label{andre2111}\emph{(\cite[th. III.2.1.11, prop. III.2.1.13]{andre1})} Assume $K$ algebraically closed. Let $\overline S$ be a manifold, and let $Z$ be a Zariski closed nowhere dense reduced analytic subset. Then the functor from tempered coverings of $\overline S$ to tempered coverings of $S=\overline S\backslash Z$ is fully faithful. If $Z$ is of codimension $\geq 2$, this functor is an equivalence of categories.
\end{prop}
We will follow the proof of~\cite[cor. X.3.4]{sga} to get the birational invariance of the tempered fundamental group of smooth and proper $K$-schemes.\\
Recall that every rational map between smooth and proper $K$-schemes is defined above the complement of a closed subset of codimension $\geq 2$.\\
Let $f:X\to Y$ is a dominant rational map between smooth and proper $K$-schemes. $f$ is defined on a Zariski open subset $U$ of $X$ (denote by $f_U$ the morphism $U\to Y$ and by $i_U$ the immersion $U\to X$) whose complement is of codimension $\geq 2$, one gets a functor from $\Covtemp(Y)$ to $\Covtemp(X)$ and one can compose it with a quasi-inverse of $\Covtemp(U)\to\Covtemp(X)$: one thus get a functor $f^*_{(U)}:\Covtemp(Y)\to\Covtemp(X)$, such that $i_U^*f^*_{(U)}$ is isomorphic to $f^*_U$. If one takes another Zariski open subset $U'$ of $X$ satisfying the same properties, one gets that $i_{U\cap U'}^*f^*_{(U)}$ and $i_{U\cap U'}^*f^*_{(U')}$ are both isomorphic to $f^*_{U\cap U'}$, and thus $f^*_{(U)}$ and $f^*_{(U')}$ are isomorphic, since $X\backslash U\cap U'$ is also of codimension $\geq 2$. Thus one gets an outer homomorphism of topological groups $f_*:\gtemp(X)\to\gtemp(Y)$, which does not depend on $U$. In particular if $f$ is a morphism of schemes, one can choose $U=X$ and thus $f_*$ is the usual outer morphism $\gtemp(X)\to\gtemp(Y)$.\\
Let $g:Y\to Z$ be another dominant rational map between smooth and proper $K$-scheme, it is defined on a Zariski open subset $V$ of $Y$ of codimension $\geq 2$, and $gf:X\to Z$ is also a dominant rational map between smooth and proper schemes, so it is also defined over a Zariski open subset $W$ of $X$ of codimension $\geq 2$. Let $U_0=U\cap f_U^{-1}(V)\cap W$ ($X\backslash U_0$ does not need to be of codimension $\geq 2$). There are morphisms $U_0\to V$ and $V\to Z$ representing $f$ and $g$ such that the composed morphism $(gf)_{U_0}:U_0\to Z$ represents $gf$. One then gets that $i_{U_0}^*f^*_{(U)}g^*_{(V)}$ and $i_{U_0}^*(gf)^*_{(W)}$ are both isomorphic to $(gf)_{U_0}^*$. Since $i_{U_0}^*$ is fully faithful, $f^*_{(U)}g^*_{(V)}$ and $(gf)^*_{(W)}$ are isomorphic (and $g_*f_*=(gf)_*$). Thus one gets a functor from the category of smooth and proper $K$-schemes with dominant rational maps to the category of groups with outer homomorphisms.\\
In particular, 
\begin{prop}\label{birat} Let $X\to Y$ be a birational map between smooth and proper $K$-schemes. Then
\[\gtemp(X)\to\gtemp(Y)\]
is an isomorphism.\end{prop}

If one considers the full subcategory $\Covtemp(X)^{(p')}$ of tempered coverings that become topological after pullback by a finite Galois covering of order prime to $p$ where $p$ is the residual characteristic, one gets in the same fashion a pro-${(p')}$ version $\gtemp(X,\bar x)$ of the tempered fundamental group (see~\cite[Rm. 3.10.1.]{mochi} in the case of a curve).

\subsection{Mochizuki's results on the pro-$(p')$ tempered group of a curve}\label{mochipart}
Following~\cite[Appendix]{mochi3}, a \emph{semigraph} $\mbb G$ will be given by a set $\mcal V$ of vertices, a set $\mcal E$ of edges and, for every $e\in\mcal E$ a set with cardinality $\leq 2$ of branches $\mcal B_e$ with a $\zeta_e:\mcal B_e\to\mcal V$ (we will say that a branch $b$ of $e$ \emph{abuts} to $v$ if $\zeta_e(b)=v$). We will say that a semigraph is a \emph{graph} if for every $e\in\mcal E$, $e$ has exactly two branches.\\
Recall from~\cite[def. 2.1]{mochi} that if one has a semigraph $\mcal G$, a structure of \emph{semigraph of anabelioid} on this graph corresponds to the following data: for each vertex or edge $x$, a Galois category (also named connected anabelioid in~\cite{mochi}) $\mcal G_x$, and for each branch $b$ of edge $e$ abutting to a vertex $v$, a morphism of anabelioids (\emph{i.e.} an exact functor) $b_*:\mcal G_e\to \mcal G_v$). Recall also that working with Galois categories is equivalent to working with profinite groups up to inner automorphism.\\
A \emph{covering} $S$ of $\mcal G$ consists of data $(S_v,\phi_e)$, where for a vertex $v$, $S_v$ is an object of $\mcal G_v^{\top}$ (the category of arbitrary disjoint unions of connected objects of $\mcal G_v$, that is the topos $\gf(\mcal G_v)-\Set$), and for every edge $e$ with branches $b_1$ and $b_2$ abutting to $v_1$ and $v_2$ $\phi_{e}$ is an isomorphism between $b_{1*}S_{v_1}$ and $b_{2*}S_{v_2}$. One has a natural notion of morphism of such coverings, so that one gets a category $\Bcov(\mcal G)$. Mochizuki associates canonically to such a covering a semigraph of anabelioids $\mcal G'$ above $\mcal G$.\\
An object of $\Bcov(\mcal G)$ is finite if each $S_v$ is in $\mcal G_v$, \emph{topological} if for each $v$ $S_v$ is a constant object of $\mcal G_v^{\top}$, and \emph{tempered} if it becomes topological after pulling back along some finite covering. The full subcategory of tempered coverings is then denoted by $\Btemp(\mcal G)$ (and $\Ba(\mcal G)$ denotes the full subcategory of finite coverings; if $\mbb G$ is connected, then it is a Galois category of fundamental group denoted as $\ga(\mcal G)$).\\
Then, if $v$ is a vertex of $\mbb G$ and $F$ is a fundamental functor of $\mcal G_v$ (and thus extends to a point of $\mcal G_v^{\top}$), one can define a functor $F_{(v,F)}:\Bcov(\mcal G)\to\Set$ which maps $S$ to $F(S_v)$ (if one changes the base point $(v,F)$, one gets an isomorphic functor), and let $F^{\temp}_{(v,F)}$ be its restriction to $\Btemp(\mcal G)$. Then we define $\gtemp(\mcal G,(v,F))$ to be $\Aut(F^{\temp}_{(v,F)})$.\\

Let $(\overline X,D)$ be a smooth $n$-pointed curve of type $g$ over a complete discrete valuation field $K$, let $X=\overline X\backslash D$, let $(\overline{\mcal X},\mcal D)$ be a semistable model over $O_{\overline K}$, where $\overline K$ is the completion of an algebraic closure of $K$, and let $\mcal X=\overline{\mcal X}\backslash\mcal D$.\\
The semigraph $\mbb G_{\mcal X}^c$ of $\mcal X_s$ is defined as follows: the vertices are the irreducible components of $\mcal X_s$, the edges are the nodes and the marked points. A node $e$ has two branches that abut to the irreducible components that contain $e$; a marked point $e$ has only one branch that abuts to the irreducible component containing the marked point.\\
When $\mcal X$ is the stable model of $X$, $\mbb G_{\mcal X}^c$ will simply be denoted by $\mbb G_{X}^c$ (or $\mbb G^c$ when there is no risk of confusion).\\
One can endow the semigraph $\mbb G^c$ of the special fiber $\mcal X_s$ of $\mcal X$ with a structure of "semigraph of anabelioids" $\mcal G^c$.
Indeed for a vertex $v_i$ corresponding to an irreducible component $C_i$ of $\mcal X_s$, consider the open subset $U_i$ of the normalization 
$\overline C'_i$ of $\overline C_i$ which is the complementary of the marked points and of the preimages of the double points of $\mcal X_s$ (the points of $\overline C'_i-U_i$ thus
correspond exactly to branches abutting to $v_i$). Then the group $\Pi_{(v_i)}$ is the tame fundamental group $\gt(U_i)$ of $U_i$ in $\overline C'_i$. The group of an edge
is $\widehat{\mathbf Z(1)}^{(p')}=\varprojlim_{(n,p)=1} \mu_n (\simeq \widehat{\mathbf Z}^{(p')})$ ($(p')$ is for the pro-prime-to-$p$ maximal quotient), 
which is canonically isomorphic to the subgroup of monodromy in $\gt(U_i)$ of a point in $\overline C'_i -U_i$. The morphism corresponding to a branch is the embedding of the monodromy group of the corresponding point of $\overline C'_i-U_i$ (defined up to conjugation), 
whereas, for an edge with two branches, one identifies the two $\widehat{\mathbf Z(1)}^{(p')}$ by $x \mapsto x^{-1}$.\\
If $\mcal G^{(p')}$ denotes the semigraph of anabelioids obtained from $\mcal G$ by replacing each profinite group by its pro-$(p')$ completion, and if $\gtemp(X_{\overline K})^{(p')}$ is the pro-${(p')}$ version of the tempered fundamental group of $X_{\overline K}$, 
\[\label{revkumm} \gtemp(\mathcal G_X^{c,(p')})=\gtemp(X_{\overline K})^{(p')}.\]
Mochizuki then shows:
\begin{thm}\emph{(\cite[cor. 3.11]{mochi})}\label{mochi} If $X_{\alpha}$ and $X_\beta$ are two curves, every isomorphism $\gamma:\gtemp(X_{\alpha,\overline K})\simeq\gtemp(X_{\beta,\overline K})$ determines, functorially in $\gamma$ up to 2-isomorphism, an isomorphism of semigraphs of anabelioids $\gamma':\mcal G^c_{X_\alpha}\simeq\mcal G^c_{X_\beta}$.\end{thm}
More precisely, the following induced diagram of topological group is commutative:
\[\xymatrix{\gtemp(\mathcal G_{X_\alpha}^{c,(p')}) \ar@{=}[d] \ar[r]^{\gtemp(\gamma'^{(p')})} & \gtemp(\mathcal G_{X_\beta}^{c,(p')})\ar@{=}[d]\\ \gtemp(X_{\alpha,\overline K})^{(p')} \ar[r]^{\gamma^{(p')}} & \gtemp(X_{\beta,\overline K})^{(p')}}\]
The vertices of the graph then correspond to the conjugacy classes of maximal compact subgroups of $\gtemp(X)^{(p')}$ (thus called \emph{verticial} subgroups of $\gtemp(X)^{(p')}$), and edges to conjugacy class of nontrivial intersections of verticial subgroups.

\section{Abelianized tempered fundamental group of a curve}
Here, we will be interested in the tempered fundamental group of an abelian
variety and in the abelianized tempered fundamental group of a curve.\\
We will first prove that if $A$ is an abelian variety over an algebraically closed complete nonarchimedean field of characteristic $0$, the tempered fundamental group is abelian and fits in a split exact sequence:
\[0\to T\to \gtemp(A)\to\gtop(A)\to 0\]
where $T$ is profinite (and more precisely isomorphic to $\widehat{\mbf Z}^n$ for some $n$).\\
Then we will prove that (as in the case of algebraic fundamental groups or complex topological fundamental groups) if $C$ is a curve and $A$ its Jacobian variety, the natural morphism $\gtemp(C)^{\ab}\to\gtemp(A)$ is an isomorphism. 
\subsection{Tempered fundamental group of an abelian variety}\label{abvar}
Let $K$ be an algebraically closed complete nonarchimedean field of
characteristic $0$.\\
Let $A$ be an abelian variety over $K$, and let $g$ be its dimension.\\

Recall the basics of $p$-adic uniformization of abelian varieties.\\
By~\cite{FvdP}, there is a commutative algebraic group $G$ (more
precisely a semiabelian variety) and a surjective analytic morphism
$u:G^{\an}\to A^{\an}$ which is the universal topological covering of $A$
and $\ker u$ is a discrete free abelian group $\Lambda$ of rank $d$.\\
Moreover, we know that $(A^{(n)}\to A)_{n\in\mathbf N}$ ($\mathbf N$ being
ordered by divisibility), with $A^{(n)}$ a copy of $A$ and $A^{(n)}\to A$
multiplication by $n$, is a cofinal family of finite Galois coverings of
$A$ (by~\cite[lecture XI]{sga}).\\
Let $G^{(n)}$ be the universal topological covering of $A^{(n)}$ (which is
isomorphic to $G$ since $A^{(n)}$ is isomorphic to $A$): one has
\[\gtemp(A)=\varprojlim_n \Gal(G^{(n)}/A).\]\\
Since $\gtop(A^{(n)})=\gtop(A)=\Lambda \simeq \mathbf Z^d$ is residually finite and is a subgroup of finite index of $\Gal(G^{(n)}/A)$ , $\Gal(G^{(n)}/A)$ is residually finite for every $n$, so $\gtemp(A)$ is also topologically residually finite as a projective limit of residually finite groups. Thus
$\gtemp(A)\to\ga(A)$ is injective and, since $\ga(A)$ is abelian, $\gtemp(A)$
is also an abelian group.\\

If $n|m$, one has the following commutative diagram :
\[\xymatrix{0 \ar[r] & \Gal(G^{(m)}/G) \ar[d] \ar[r] & \Gal(G^{(m)}/A)
  \ar[d] \ar[r] & \Gal(G/A)=\Lambda \ar@{=}[d] \ar[r] & 0\\
0 \ar[r] & \Gal(G^{(n)}/G) \ar[r] & \Gal(G^{(n)}/A) \ar[r] &
\Gal(G/A)=\Lambda \ar[r] & 0}.\]

Let us write $T(G)=\varprojlim \Gal(G^{(n)}/G)$. This is a profinite abelian
group, so it splits (canonically) as the product of its pro-$p$-Sylow
subgroups $T_l(G)$.\\
By taking the projective limit in the previous commutative diagram, one gets the following axact sequence
(it is right exact because the $\Gal(G^{(n)}/G$ are finite): 
\[\xymatrix{0 \ar[r] & T(G) \ar[r] & \gtemp(A)
  \ar[r] & \Lambda \ar[r] & 0},\label{seA}\]
with $\gtemp(A)$ abelian. Thus it is an exact sequence of abelian groups.\\
But $\Lambda$ is a free abelian group, so the exact sequence must be
split. One thus gets a non canonical isomorphism:
\[\gtemp(A)\simeq \Lambda \times T(G).\]

By taking the pro-$l$ completion of the isomorphism here above, one gets: 
\[\mathbf Z_l^{2g}\simeq\ga(A)^l\simeq\widehat{\gtemp(A)}^l\simeq\mathbf Z_l^d\times T_l(G),\]
and so $T_l(G)\simeq\mathbf Z_l^{2g-d}$.
We finally obtained a non canonical isomorphism:
\[\gtemp(A)\simeq\mathbf Z^d\times\widehat{\mathbf Z}^{2g-d}.\]

\subsection{Jacobian variety and $\gtemp^{\ab}$ of a curve}
If $G$ is a topological prodiscrete group with a countable basis of
neighborhoods of $1$, $G^{\ab}$ will be the
topological group $G/\overline{D(G)}$ (where $\overline D(G)$ is the
closure of the derived subgroup of $G$), it is also a prodiscrete group and
$G\to G^{\ab}$ makes $G^{\ab}-\Ens$ as a full subcategory of $G-\Ens$.\\

Let $K$ be a complete discrete valuation field of characteristic $0$, let
$\bar K$ be the completion of its algebraic closure, let $C$ be a curve over $K$
and let $A$ be the Jacobian variety of $C_{\bar K}$.\\
Let $P$ be a closed point of $C_{\bar K}$.\\
Consider the morphism $C_{\bar K}\to A$ that maps $x$ to the divisor $[x]-[P]$. One gets a homomorphism $\gtemp(C_{\bar K},P)^{\ab}\to\gtemp(A,0)$ that factorizes through $\gtemp(C_{\bar K},P)^{\ab}$ since $\gtemp(A,0)$ is abelian.\\
\begin{thm}\label{jacobian} The morphism $\gtemp(C_{\bar K},P)^{\ab}\to\gtemp(A,0)$ is an isomorphism.\end{thm}
\dem
We have a morphism $C_{\bar K}^g\to A$ which maps $(x_1,\dots,x_g)$ to the
divisor $[x_1]+\cdots+[x_g]-g[P]$ of $C_{\bar K}$ . This morphism is
invariant by the action of $\frak S_g$ on $C_{\bar K}^g$ and thus factorizes in a morphism $C_{\bar K}^{(g)}:=C_{\bar K}^g/\frak S_g\to
A$. Remind that this is a birational morphism and that $C_{\bar K}^{(g)}$
is smooth over $\bar K$ (voir~\cite[Th. 5.1.(a), prop. 3.2]{jacmilne}).\\
We thus get a sequence of morphisms:
\[C_{\bar K}\to C_{\bar K}^g\to C_{\bar K}^{(g)}\to A,\]
where the left morphism maps $x$ to $(x,P,\dots, P)$ and where the 
composed morphism is the morphism that maps $x$ to $[x]-[P]$.\\

Since $C_{\bar K}^{(g)}\to A$ is a birational morphism of proper smooth $\bar
K$-varieties, $\gtemp(C_{\bar
  K}^{(g)},(P,\dots,P))\to\gtemp(A,0)$ is an isomorphism according to proposition~(\ref{birat}).\\
Thus $\gtemp(C_{\bar K}^{(g)},(P\dots,P))$ is abelian and residually finite.\\
Thus $\gtemp(C_{\bar K},P)\to\gtemp(C_{\bar K}^{(g)},(P,\dots,P))$
factorizes through $\phi:\gtemp(C_{\bar K},P)^{\ab}\to\gtemp(C_{\bar
  K}^{(g)},(P,\dots,P))$.\\

As $\gtemp(C_{\bar K},P)^{\ab}$ is a projective limit of abelian groups of finite type (and thus residually finite groups), it must be residually finite too.\\
The following commutative diagram:
\[\xymatrix{\gtemp(C_{\bar K},P)^{\ab} \ar[r]^{\phi} \ar[d] & \gtemp(C_{\bar
    K}^{(g)},(P,\dots,P)) \ar[d]\\ \ga(C_{\bar K},P)^{\ab} \ar[r]^{\simeq} & \ga(C_{\bar
    K}^{(g)},(P,\dots,P))},\]
whose vertical arrows are injective, shows that $\phi$ is
injective.\\

One may also put an orbifold structure on $C_{\bar K}^{(g)}$, as being $C_{\bar K}^g/\frak S_g$ (then $\gorb(C^{(g)}_{\overline K})$ denotes the corresponding extension of $\frak S_g$ by $\gtemp(C^g_{\overline K})$, as in~\cite[prop. III.4.5.8]{andre1}).\\
One then has the following commutative diagram whose line is exact:
\[\xymatrix{
1 \ar[r] & \gtemp(C_{\bar K}^g,(P,\dots,P)) \ar[r]^i \ar[ddr]_\alpha &
\gorb(C_{\bar K}^{(g)},(P,\dots,P)) \ar[r]\ar[d]^{\pi_1} & \frak S_g \ar[r] & 1\\
& & \gorb(C_{\bar K}^{(g)},(P,\dots,P))^{\ab} \ar[d]^{\pi_2} 
& & \\  & & \gtemp(C_{\bar K}^{(g)},(P,\dots,P)) & &}\]
The morphisms $i$, $\pi_1$ and $\pi_2$ are open, thus $\alpha=\pi_2 \pi_1 i$ is open: $\Image \alpha$ is an open subgroup of $\gtemp(C_{\bar K}^{(g)},(P,\dots,P))$.\\
But, if it is a strict subgroup, $\gtemp(C_{\bar K}^{(g)},(P,\dots,P))/\Image \alpha$ is a nontrivial abelian group of finite type and thus has a nontrivial finite quotient, which corresponds to an open subgroup of finite index of $\gtemp(C_{\bar K}^{(g)},(P,\dots,P))$ which contains $\Image \alpha$.\\
But $\hat \alpha$ is surjective (since $\ga(C_{\bar K})\to\ga(C_{\bar K}^{(g)})$, which factorizes through $\hat \alpha$ is surjective), and thus one gets a contradiction: $\alpha$ is also surjective.\\

Let us now consider the diagonal hommomorphism $\delta:\gtemp(C_{\bar K},P)\to\gtemp(C_{\bar
K}^{g},(P,\dots,P))$. By identifying $\gtemp(C_{\bar
K}^{g},(P,\dots,P))$ to $\gtemp(C_{\bar
K},(P,\dots,P))^g$ thanks to proposition~(\ref{prod}) that we will prove in the next section, $\delta$ may be identified to $\gtemp(C_{\bar
K},P)\to\gtemp(C_{\bar K},P)^g:g\mapsto (g,1,\dots,1)$.\\
Thus $\gtemp(C_{\bar
K}^{g},(P,\dots,P))$ is generated by the images of $\sigma\circ\delta$
when $\sigma$ describes $\frak S_g$, and as $\alpha$ is invariant by
$\sigma$, $\alpha\circ\delta$ is surjective, so $\gtemp(C_{\bar
  K},P)^{\ab}\to\gtemp(C_{\bar K}^{(g)},(P,\dots,P))$ is also surjective,
thus it is bijective.\\
If $U$ is an open subgroup of $\gtemp(C_{\bar K},P)$, the
group generated by the $\sigma(\delta(U))$ is an open subgroup of $\gtemp(C_{\bar
K}^{g},(P,\dots,P))$, thus, as $\alpha$ is open and $\frak
S_g$-invariant, $\alpha\circ\delta$ is open, thus $\gtemp(C_{\bar
  K},P)^{\ab}\to\gtemp(C_{\bar K}^{(g)},(P,\dots,P))$ is also open, so it is an isomorphism.\findem

\section{Alterations and tempered fundamental group}
In this section we will describe two applications of the theorems of
existence of semistable alterations proved by de Jong to the tempered
fundamental group. Indeed, Berkovich already showed in~\cite[§9]{berk2} how
those alterations help building a skeleton, homeomorphic to the geometric
realization of a polysimplicial set, on which a Zariski open set of the
variety retracts.\\
We will here deduce from the results of Berkovich that, in characteristic
zero, the tempered fundamental group of a smooth algebraic variety is
invariant by base change of algebraically closed complete fields, and that
the tempered fundamental group of the product of two smooth varieties (over
an algebraically closed base) is canonically isomorphic to the product of
the tempered fundamental groups of each variety.
\subsection{Preliminaries about the skeleton of a Berkovich space}\label{berkspaces}
Let $k$ be a complete nonarchimedean field and let $k^{\circ}$ be its ring of integers.\\
Recall the definition of a  polystable morphism of formal schemes:
\begin{dfn}\emph{(\cite[def. 1.2]{berk2},\cite[section 4.1]{berk3})} Let $\phi:\fk X\to \fk Y$ be a locally finitely presented morphism of formal schemes.
\begin{enumerate}[(i)]
\item $\phi$ is said to be \emph{strictly polystable} if, for every point $y\in\fk Y$, there exists an open affine neighborhood $\fk X'=\Spf(A)$ of $\phi(y)$ and an open neighborhood $\fk Y'\subset\phi^{-1}(\fk X')$ of $x$ such that the induced morphism $\fk Y'\to \fk X'$ goes through an étale morphism $\fk Y'\to\Spf(B_0)\times_{\fk X'}\cdots\times_{\fk X'}\Spf(B_p)$ where each $B_i$ is of the form $A\{T_0,\cdots,T_{n_i}\}/(T_0\cdots T_n-{a_i})$ with $a\in A$ and $n\geq 0$. It is said to be \emph{nondegenerate} if one can choose $X'$, $Y'$ and $(B_i,a_i)$ such that $\{x\in(\Spf(A)_\eta)|a_i(x)=0\}$ is nowhere dense.
\item $\phi$ is said to be \emph{polystable} if there exists a surjective étale morphism $\fk Y'\to\fk Y$ such that $\fk Y'\to \fk X$ is strictly polystable. It is said to be \emph{nondegenerate} if one can choose $\fk Y'$ such that $\fk Y'\to \fk X$ is nondegenerate. \end{enumerate}\end{dfn}
Then a (\emph{nondegenerate}) \emph{polystable fibration} of length $l$ over $\fk S$ is a sequence of (nondegenerate) polystable morphisms $\underline{\fk X}=(\fk X_l\to\cdots\to\fk X_1\to\fk S)$.\\

Berkovich defines \emph{polysimplicial sets} in~\cite[section 3]{berk2} as follows.\\
For an integer $n$, denote $[n]=\{0,1,\cdots,n\}$.\\
For a tuple $\mbf n=(n_0,\cdots,n_p)$ with either $p=n_0=0$ are $n_i\geq 1$ for all $i$, let $[n]$ denote the set $[n_0]\times\cdots\times [n_p]$ and $w(\mbf n)$ will denote the number $p$.\\
Berkovich defines a category $\mbf \Lambda$ whose objects are $[\mbf n]$ and morphisms are maps $[\mbf m]\to [\mbf n]$ associated with triples as follows. $J$ is a subset of $[w(\mbf m)]$ assumed to be empty if $[\mbf m]=[0]$, $f$ is an injective map $J\to [w(\mbf n)]$ and $\alpha=\{\alpha_l\}_{0\leq l\leq p}$, where $\alpha_l$ is an injective map $[m_{f^{-1}(l)}]\to [n_l]$ if $l\in\Image(f)$, and $\alpha_l$ is a map $[0]\to [n_l]$ otherwise. The map $\gamma:[\mbf m]\to [\mbf n]$ associated with $(J,f,\alpha)$ takes $\mbf j=(j_0,\cdots,j_{w(m)})\in [m]$ to $\mbf i=(i_0,\cdots,i_{w(\mbf n)})$ with $i_l=\alpha_l(j_{f^{-1}(l)})$ for $l\in\Image(f)$, and $i_l=\alpha_l(0)$ otherwise.\\
Then a polysimplicial set is a functor $\mbf \Lambda^{\op}\to\Set$, they form a category denoted $\mbf \Lambda^{\circ}\Set$.\\
Let $\Delta$ be the strict simplicial category of integers with increasing maps.
Strict simplicial sets are simple examples of polysimplicial sets by extending to $\Delta^{\circ}\Set\to\mbf\Lambda^{\circ}\Set$ the functor $\Delta\to \mbf \Lambda$, that maps $[n]$ to itself, such that it commutes with direct limits.\\
Berkovich then considers a functor $\Sigma:\mbf \Lambda\to\mcal Ke$ (where $\mcal Ke$ is the category of Kelley spaces, \emph{i.e.} topological spaces $X$ such that a subset of $X$ is closed whenever its intersection with any compact subset of $X$ is closed) that takes $[\mbf n]$ to $\Sigma_{\mbf n}=\{(u_{il})_{0\leq i\leq p,0\leq l\leq n_i}\in [0,1]^{[\mbf n]}|\sum_l u_{il}=1\}$, and takes a map $\gamma$ associated to $(J,f,\alpha)$ to $\Sigma(\gamma)$ that maps $\mbf u=(u_{jk})$ to $\mbf u'=(u'_{il})$ defined as follows: if $[\mbf m]\neq [0]$ and $i\notin\Image(f)$ or $[\mbf m]=[0]$ then $u'_{il}=1$ for $l=\alpha_i(0)$ and $u'_{il}=0$ otherwise; if $[\mbf m]\neq [0]$ and $i\in\Image(f)$, then then $u'_{il}=u_{f^{-1}(i),\alpha_i^{-1}(l)}$ for $l\in\Image(\alpha_i)$ and $u'_{il}=0$ otherwise. Berkovich thus gets a functor, the \emph{geometric realization}, $|\ |:\mbf \Lambda^{\circ}\Set\to\mcal Ke$ by extending $\Sigma$ such that it commutes with direct limits.\\

Berkovich attached to a polystable fibration over $\underline{\fk X}=(\fk X_l\to\fk X_{l-1}\to\cdots\to\Spf(k^{\circ}))$ a polysimplicial set $\mbf C_l(\fk X_s)$ (which only depends on the special fiber) and a subset of the generic fiber $\fk X_{l,\eta}$ (this is a Berkovich space) of $\fk X_l$ named the skeleton of $\underline{\fk X}$ and denoted by $S(\underline{\fk X})$ which is canonically homeomorphic to $|\mbf C_l(\fk X_s)|$ (see~\cite[th. 8.2]{berk2}), such that $\fk X_{l,\eta}$ retracts by a proper strong deformation onto $S(\underline{\fk X})$.\\

In fact, when $\underline{\fk X}$ is non degenerate (for example generically smooth, we will in fact only use the results of Berkovich to such polystable fibrations), the skeleton of $\underline{\fk X}$ only depends on $\fk X_l$ (such a formal scheme that fits into a polystable fibration will be called a \emph{pluristable morphism}, and we will note $S(\fk X_l)$ this skeleton) according to~\cite[prop. 4.3.1.(ii)]{berk3}. In this case~\cite[prop. 4.3.1.(ii)]{berk3} give a simple explicit description of $S(\fk X_l)$. For any $x,y\in\fk X_{l,\eta}$, we note $x\preceq y$ if for every $\fk X'\to\fk X_l$ and every point $x'$ over $x$, there exists $y'$ over $y$ such that for every $f\in O(\fk X_\eta)$, $|f(x')|\leq |f(y')|$ ($\preceq$ is a partial order on $\fk X_{l,\eta}$). Then $S(\fk X_l)$ is just the set of maximal points of $\fk X_{l,\eta}$ for $\preceq$.\\

We will not give in full detail the construction of $\mbf C_l(\fk X_s)$ (we will not need it).
Rather, we will give some examples for a polystable fibration of length 1 (they are all easily deduced from the construction explained in~\cite[section 3]{berk2}).\\
If $\fk X$ is just $\Spf(B_0)\times_{\Spf k^{\circ}}\cdots\times_{\Spf k^{\circ}}\Spf(B_{p+1})$ where each $B_i$ is of the form $k^{\circ}\{T_0,\cdots,T_{n_i}\}/(T_0\cdots T_{n_i}-a_i)$ with $a_i\in k^\circ$, $|a_i|<1$ for $i<p$ and $|a_{p+1}|=1$, then $\mbf C_1(\fk X\to\Spf(k^{\circ}))$ is just $[\mbf n]$ with $\mbf n=(n_0,\cdots,\mbf n_p)$.\\
If $\fk X$ is just a semistable formal scheme over $\Spf(k^{\circ})$, then $\mbf C_1(\fk X\to\Spf(k^{\circ}))$ is a  simplicial set.\\
For example, if $\mcal X$ is a semistable model of a complete curve (then it is clearly polystable over $k^{\circ}$), the polysimplicial set defined by Berkovich is just the graph of the stable reduction defined in the previous section (a graph can be considered as a simplicial set of dimension 1, and thus as a polysimplicial set).\\
If $(\overline{\mcal X},\mcal D)$ is a semistable model of a curve (not necessarily complete), the graph of $(\overline{\mcal X})$ is the semigraph of $\mcal X=\overline{\mcal X}\backslash\mcal D$ in which one deleted all the edges with only one branch. The retraction of $(\overline{\mcal X}_\eta)^{\an}$ to the geometric realization of this graph restricts to a retraction of $({\mcal X}_\eta)^{\an}$ thanks to~\cite[cor. 8.4]{berk2}.
\\

The retraction to $S(\underline{\fk X})$ commutes with étale morphisms:
\begin{thm}\label{berk81}\emph{(\cite[th. 8.1]{berk2})}
One can construct for every polystable deformation $\fk X=(\fk
X_l\stackrel{f_{l-1}}{\to}\cdots\stackrel{f_1}{\to}\fk X_1\to\Spf(k^{\circ}))$ a proper strong deformation retraction $\Phi^l:\fk X_{l,\eta}\times [0,l]\to \fk X_{l,\eta}$ of $\fk
X_{l,\eta}$ onto the skeleton $S(\underline{\fk X})$ of $\underline{\fk X}$ such that :\begin{enumerate}[(i)]
\item $S(\underline{\fk X})=\bigcup_{x\in S(\underline{\fk X}_{l-1})}S(\fk X_{l,x})$, where $\underline{\fk X}_{l-1}=(\fk X_{l-1}\to\cdots\to\Spf(k^{\circ}))$, and
  $f_{l-1,\eta}(x_t)=f_{l-1,\eta}(x)_{t-1}$ for every $t\in [1,l]$;
\item $(x_t)_{t'}=x_{\max(t,t')}$;
\item $x\leqslant x_t$;
\item for every $x$ and $i\leqslant l-1$, there is $t'\in [i,i+1]$ such that
  $x_t=x_i$ if $t\in [i,t']$ and $t\mapsto x_t$ is injective over $[t',i+1]$;
\item for every $t$, $\pi(x_t)$ is contained in the same stratum of $\fk
  X_{l,s}$ as $\pi(x)$; moreover $\pi(x_l)$ is the generic point of this stratum
  and $\pi(x_i)=\pi(x_t)$ for every $0\leqslant i\leqslant l-1$
  and every $t\in [i,i+1]$;
\item if $\phi: \fk Y\to \fk X$ is a morphism of fibrations in $\mcal
  Pst f^{\et}_l$, one has $\phi_{l,\eta}(y_t)=\phi_{l,\eta}(y)_t$ for every
  $y\in\fk Y_{l,\eta}$.\end{enumerate}\end{thm}
Berkovich deduces from the functoriality of the retraction with respect to \'etale morphisms the following corollary:
\begin{cor}\label{berk85}\emph{(\cite[cor 8.5]{berk2})} Let $\underline{\fk X}$ be a polystable fibration over   
  $k'^\circ$ (where $k'$ is a finite normal extension of
  $k$) with a normal generic fiber $\fk X_{l,\eta}$. Suppose we are given an action
  of a finite group $G$ on $\underline{\fk X}$ over $k^\circ$ and a Zariski open dense subset $U$ of $\fk X_{l,\eta}$. 
  Then there is a strong deformation retraction of the Berkovich space
  $G\backslash U$ to a closed subset homeomorphic to $G\backslash |\mbf C^l(\underline{\fk X})|$.\end{cor}
More precisely, in this corollary, the closed subset in question is the image of $S(\underline{\fk X})$ (which is $G$-equivariant and contained in $U$) by $U\to G\backslash U$.\\

Recall that the skeleton is functorial for pluristable morphisms:
\begin{prop}\emph{(\cite[prop. 4.3.2.(i)]{berk3})} If $\phi:\fk X\to\fk Y$ is a pluristable morphism between nondegenerate pluristable formal schemes over $k^\circ$, $\phi_\eta(S(\fk X))\subset S(\fk Y)$.\end{prop}
In fact, more precisely, from the construction of $S$, $S(\fk Y)=\bigcup_{x\in S(\fk X)}S(\fk Y_x)$.
Recall also that, if $\bar k$ is algebraically closed, the polysimplicial complex of a polystable fibration commutes with base change:
\begin{prop}\label{berk610}\emph{(\cite[prop. 6.10]{berk2})} If $\underline X$ is a polystable fibration over a field $K$, there exists a finite separable extension $K'$ of $K$ such that for any bigger field $K''$, $\mbf C(\underline X\otimes K'')\to\mbf C(\underline X\otimes K')$ is an isomorphism.\end{prop}

In order to use the previous description of the Berkovich space of a scheme with a model over $k^{\circ}$ which admits a polystable fibration for understanding the topology of a smooth scheme over $k$, we will need de Jong's result about the existence of alterations by such pluristable schemes over $k^{\circ}$.\\
More precisely we will use the following consequence of de Jong's theory given by Berkovich (as we will work over valuation fields of characteristic 0, we give here only a restricted version in this case):
\begin{lem}\label{berk92} \emph{(\cite[lem. 9.2]{berk2})} Let $A$ be a Henselian valuation ring with fraction field of characteristic 0, let $X$ be an integral  scheme proper flat and finite presentation over
  $A$, with an irreducible generic fiber of dimension $l$. Then there is:\begin{enumerate}[(a)]
\item a finite Galois extension of the fraction field of $A$, with ring of integers $A'$,
\item a polystable fibration $\underline X'=(X'_l\to\cdots\to X'_0=\Spec
  A')$, where every morphism is projective of dimension 1 with smooth
  geometrically irreducible generic fibres,
\item an action of a finite group $G$ on $\underline X'$ over $A$,
\item a commutative diagram
\[\xymatrix{X'_l\ar[r]^\phi\ar[d] & X\ar[d] \\ \Spec(A')\ar[r] &
  \Spec(A)}\]
where $\phi$ is a dominant $G$-equivariant morphism such that the generic fiber is generically étale with Galois  group $G$.\end{enumerate}\end{lem}

\subsection{Invariance of $\gtemp$ by base change of algebraically closed
  fields}
Let $X$ be a smooth algebraic variety over an algebraically closed complete
nonarchimedean field $K$ of characteristic 0.\\
Let $X_i\to X$ be a finite étale connected Galois covering and $U_0$ a
dense affine Zariski open subset of $X_i$, and let us embed ({\it{e.g.}}
by~\cite[lem 9.4]{berk2}) $U_0$ in a scheme $\overline{X_i}$ which is
proper , of finite presentation and flat over $O_K$.\\
Then there is by lemma~(\ref{berk92}) a generically smooth polystable
fibration $X'_i$ over $O_K$ endowed with an action of a group $G$ such that
$(X'_i,G)$ is a Galois alteration of $\overline{X_i}$.\\
Let $U$ be a dense Zariski open subset of $\overline{X_i}$ included in
$U_0$ (and therefore in $X_i$) such that $U'\to U$ is finite (where $U'$ is
the pullback of $U$ in $X'_i$).\\
Then $U^{\an}$ retracts by strong deformation on $G\backslash S(X'_{i,s})$ by corollary~(\ref{berk85}).\\

Let $K'/K$ an algebraically closed extension of complete valued fields and let $x'$ a geometric point with image.\\
$X_{i,K'}\to X_{K'}$ is also a finite étale connected Galois covering
(by~\cite[cor X.1.8]{sga}), $X'_{i,0_{K'}}$ is also a polystable fibration
endowed with an action of $G$ and $(X'_{i,O_{K'}}$ is also a Galois
alteration of $\overline{X_{i}}_{,K'}$, finite over $U_{K'}$.\\
Thus, as in the previous case, $U^{\an}_{K'}$ retracts on the closed subset
 $G\backslash S(X'_{i,K'})$ (and the natural
morphism $U^{\an}_{K'}\to U^{\an}$ maps $S(X'_{i,K',s})$ to
$S(X'_{i,s})$.\\
But $\mathbf C(X'_{i,K',s})\to\mathbf C(X'_{i,s})$ is an isomorphism
according to proposition~(\ref{berk610}). The morphism $U^{\an}_{K'}\to
U^{\an}$ is thus a homotopy equivalence.\\
One has the following 2-commutative diagram :
\[\xymatrix{\Covtop(U_{K'}) & \Covtop(U) \ar[l] \\ \Covtop(X_{i,K'})\ar[u]
  & \Covtop(X_i)\ar[u]  \ar[l]}.\]
The vertical arrows are equivalences by proposition~(\ref{andre114}), and
  the top arrow is also an equivalence by what has just been shown.\\
Thus, if $x'_i$ is a geometric point of $X_{i,K'}$ with image $x_i$ in $X_i$, $\gtop(X_{i,K'},x'_i)\to \gtop(X_i,x_i)$ is an isomorphism.\\
If $X^{\infty}_i$ is the universal topological covering of $X_i$ and
  $X^{\infty}_{i,K'}$ is the universal topological covering of $X_{i,K'}$,
  one has $X^{\infty}_{i,K'}=(X^{\infty}_i)_{K'}$.\\
Therefore $\Gal(X^{\infty}_{i,K'}/X_{K'})=\Gal(X^{\infty}_i/X)$.\\
By taking the projective limit over a cofinal projective system of geometrically pointed (the points being defined in a big enough valuation field $\Omega$) Galois coverings
  $((X_i)_i,x_i)_{i\in \mbf N}$ ($(X_{i,K'},x'_i)_{i\in \mbf N}$, where $x'_i$ is some point of $X_{i,K'}$ over $x_i$, is then also a cofinal projective system of
  Galois covering of $X_{K'}$ by~\cite[lecture XIII]{sga}), one gets :
\begin{prop}\label{invariance} The morphism $\gtemp(X_{K'},x')\to \gtemp(X,x)$ is an
  isomorphism. \end{prop}

\subsection{Products and tempered fundamental group}
Let $Y$ be another smooth algebraic variety over $K$, $Y_j\to Y$ a finite
étale connected Galois covering, $\overline{Y_j}$ a scheme which is proper,
of finite presentation and flat over $O_K$, in which is embedded a dense
affine open subset $V_0$ of $Y_j$, and let $(Y'_j,H)\to \overline{Y_j}$ a
Galois alteration and $Y'_j\to O_K$ a polystable fibration.\\
Let also $V\subset Y_j$ be a dense Zariski open subset, included in $V_0$
such that $V'\to V$ is finite.\\
Then, as for $X$, $V^{\an}$ retracts by strong deformation over
$H\backslash |\mathbf C(Y'_{j,s})|$.\\

One get the same results for $X_i\times Y_j$ :\\
$X'_i\times Y'_j\to X_i\times Y_j$ is a Galois alteration of group $G\times
H$, and $U\times V$ retracts on $(G\times H)\backslash S(X_i\times Y_j)$.\\
But the pluristable maps $X_i\times Y_j\to X_i$ and $X_i\times Y_j\to Y_j$ map $S(X_i\times Y_j)$ into $S(X_i)$ and $S(X_j)$ respectively. Whence a map $f:S(X_i\times Y_j)\to S(X_i)\times S(Y_j)$ (and it is compatible with the action of $G\times H$). But as 
\[S(X_i\times Y_j)=\bigcup_{x\in S(X_i)}S((X_i\times Y_j)_x)=\bigcup_{x\in S(X_i)}S(Y_j\otimes \mcal H(x))\]
and $S(Y_j\otimes\mcal H(x))\to S(Y_j)$ is an isomorphism according to proposition~(\ref{berk610}), $f$ is bijective. But $S(X_i\times Y_j)$ is compact, and thus $f$ is an isomorphism.\\
Thus $(G\times H)\backslash S(X_i\times Y_j)\to G\backslash S(X_i)\times H\backslash S(Y_j)$ is an isomorphism.\\
Thus $(U\times V)^{\an}\to U^{\an}\times V^{\an}$ is a homotopy equivalence
(the product on the right is the product of topological spaces), and
by applying proposition~(\ref{andre114}) to $U\subset X_i$, $V\subset
Y_j$ and $U\times V\subset X_i\times Y_j$, one gets that $$\gtop(X_i\times
Y_j,(x_i,y_j))\to \gtop(X_i,x_i)\times \gtop(Y_j,y_j)$$ is an isomorphism.\\
So $(X_i\times Y_j)^{\infty}=X^{\infty}_i\times Y^{\infty}_j$ and
\[\Gal((X_i\times Y_j)^{\infty}/(X\times Y))=\Gal(X_i^{\infty}/X)\times
\Gal(Y_j^{\infty}/Y).\]
If $(X_i,x_i)_i$ and $(Y_j,y_j)_j$ are cofinal projective systems of connected
geometrically pointed Galois coverings of $X$ and $Y$, $(X_i\times Y_j,(x_i,y_j))_{(i,j)}$ is a cofinal
projective system of connected Galois coverings of $X\times Y$
by~\cite[lecture XIII]{sga}.\\
By taking the projective limit over $(i,j)$ in the previous isomorphism, one
gets:
\begin{prop}\label{prod} $\gtemp(X\times Y,(x,y))\to\gtemp(X,x)\times \gtemp(Y,y)$ is an isomorphism.\end{prop}

\section{Metric structure of the graph of the stable model of a curve}
The main goal of this section is to prove that the metric structure of the graph of the stable model  (or equivalently, the graph structure of the skeleton of the curve considered as a Berkovich space, as recalled in section~\ref{berkspaces}) of a Mumford curve can be recovered from the tempered fundamental group of this curve. We will only consider the case of mixed characteristic (in the case of equal characteristic 0, this is surely false, as the whole tempered fundamental group can be recovered only from the graph of groups associated with the stable reduction of the curve according to~\cite[ex. 3.10]{mochi}).
\begin{dfn} A \emph{metric structure} over a graph $\mbb G$ is a function \[d:\{\text{edges of } \mbb G\}\to \mbf R^{>0}. \]For any $e$, $d(e)$ is called the \emph{length} of the edge $e$ with respect to the metric structure $d$.\\
A graph endowed with a metric structure is called a \emph{metric graph}.\end{dfn}
For a curve $X$ over an algebraically closed complete nonarchimedean field $k$ with a semistable model $\mcal X$, let $e$ be an edge of the graph of this semistable reduction (that is a node of $\mcal X_s$). Then etale locally at this node, $\mcal X$ is étale over $k^{\circ}[X_0,X_1]/(X_0X_1-a)$ with $a\in k^{\circ}$. According to~\cite[cor. 2.2.18]{thuillier}, $|a|$ does not depend of any choice.  Then we will denote $d(e)=-\log_p(|a|)$, which defines a natural metric structure on the graph of the stable reduction of $X$.\\
For example, if $\mcal X$ is the stable model of $\mbf P^1\backslash\{0,1,\infty,\lambda\}$ with $|\lambda|<1$ then the graph of $\mcal X$ has a single edge of length $-\log_p(|\lambda|)$.\\ 

We already know, thanks to~\cite[ex. 3.10]{mochi} that one can recover the graph of the stable model from the tempered fundamental group. In fact we will deduce from Mochizuki's study that one can recover, for every finite index open subgroup of the tempered fundamental group and every vertex of the skeleton of the curve, whether the corresponding covering of the curve is split over this vertex (an étale covering $X'\to X$ of manifolds will be said to be \emph{split} over a point $x\in X$ if, for every $x'\in X'$ over $x$, $\mcal H(x)\to\mcal H(x')$ is an isomorphism; if $X'\to X$ is a covering of order $n$, then it is split over $x$ if and only if the fiber of $x$ has cardinality $n$, if and only if, locally in a neighborhood of $x$, $X'\to X$ pulls back to a topological covering thanks to~\cite[III.1.2.1]{andre1}).\\
This suggests to consider how finite étale coverings of a Mumford curve are split over vertex points. No one knows, as far as I know, how to determine in general if a covering is split over some point of the Berkovich space. However, studying simple coverings may be enough to see that the metric structure of the skeleton must play a role in the structure of the tempered fundamental group.\\
More precisely, 
\begin{lem}\label{decompositionpuissance}The covering $\mbf G_m\stackrel{z\mapsto z^{p^e}}{\to}\mbf G_m$ is split over a Berkovich point $B(1,r)$ corresponding to the ball of center 1 and ray $r$ with $r<1$, if and only if $r<p^{-e-\frac{1}{p-1}}$.\\
More precisely, $B(1,r)$ has $p^i$ preimage if:
\begin{itemize}
\item $r\in ]p^{-\frac{p}{p-1}},1]$, if $i=0$;
\item $r\in ]p^{-i-\frac{p}{p-1}},p^{-i-\frac{1}{p-1}}]$, if $1\leq i\leq e-1$;
\item $r\in [0,p^{-e-\frac{1}{p-1}}]$, if $i=e$.
\end{itemize}\end{lem}
\dem Let $g:z\mapsto z^p$, and let us calculate $g(B(h',r))$ with $|h'|=1$ and $r<1$.
Let $f_h':z\mapsto (z+h')^p-h'^p=\sum_{i=1}^p a_iz^i$ with $|a_p|=1$ and $|a_i|=p^{-1}$ if $1\leq i\leq p-1$.\\
Then $g(B(h',r))=B(h'^p,r')$ with 
\[r'=\max_{|z|<r}|f(z)|=\max |a_i|r^i=\max\{r^p,rp^{-1}\}=\left\{\begin{array}{ll} rp^{-1} & \text{if }r\leq p^{-\frac{1}{p-1}}\\ r^p& \text{if }r\geq p^{-\frac{1}{p-1}}\end{array}\right.\]
Moreover let $h$ be of norm $1$, let $h^{1/p}$ be a $p$th root of $h$, and let $r'<1$. Then
$|h'^p-h|=\prod_{\zeta\in\mu_p}|h'-\zeta h^{1/p}|\leq r'$ implies that there exists $\zeta_0\in\mu_p$ such that $|h'-\zeta_0 h^{1/p}|\leq r'^{1/p}$ (\emph{i.e.} $B(h',r'^{1/p})=B(h^{1/p},r'^{1/p})$).\\
Suppose now $|\zeta-\zeta'|=p^{-\frac{1}{p-1}}$. Since $|\zeta-\zeta_0|=p^{-\frac{1}{p-1}}$ if $\zeta\in\mu_p\backslash\{\zeta_0\}$, $|h'-h^{1/p}\zeta|=p^{-\frac{1}{p-1}}$ if $\zeta\in\mu_p\backslash\{\zeta_0\}$ and thus $|h'-h^{1/p}\zeta_0|=|h'^p-h|/\prod_{\zeta\in\mu_p\backslash\{\zeta_0\}}|h'-h^{1/p}\zeta|\leq pr'$ (\emph{i.e.} $B(h',pr')=B(\zeta_0 h^{1/p},pr')$). Thus
\[g^{-1}(B(h,r'))=\left\{\begin{array}{ll} \{B(\zeta h^{1/p},pr')\}_{\zeta\in\mu_p} & \text{if }r'\leq p^{-\frac{p}{p-1}}\\ \{B(\zeta h^{1/p},r'^{1/p})\}_{\zeta\in\mu_p} & \text{if }r'\geq p^{-\frac{p}{p-1}}\end{array}\right.\]
Since $|\zeta-\zeta'|=p^{-\frac{1}{p-1}}$ if $\zeta\neq\zeta'\in\mu_p$, one gets that $g^{-1}(B(h,r'))$ has a single element if $r'\geq p^{-\frac{p}{p-1}}$ and $p$ otherwise. Thus ones gets the wanted result when $e=1$.\\
In the general case, one uses induction on $e$ by decomposing $z\mapsto z^{p^e}$ in $z\mapsto z^{p^{e-1}}\mapsto z^{p^e}$.\findem

Thanks to this, we will first study the case of the projective line minus some points and, by cutting and pasting this kind of covering, the case of a punctured elliptic curve.\\
For the more general case of a Mumford curve $X$, we will also study the structure of $\mbf Z/p^e\mbf Z$-torsors. The theory of theta functions as can be found in~\cite{vdp1} and~\cite{vdp2} tells us that the pullback to the universal topological covering $\Omega$ of such a torsor is in fact the pullback of $\mbf G_m\stackrel{z\mapsto z^{p^e}}{\to}\mbf G_m$ along some theta function $\Omega\to\mbf G_m$, which in turn corresponds to some equivariant current over the tree $\mbb T(\Omega)$ of $\Omega$. Therefore, we will begin our study by proving that if two currents coincide over a sufficiently big part of $\mbb T(\Omega)$, the quotient of the two corresponding invertible functions is nearly constant over some smaller part of $\mbb T(\Omega)$ and thus the two corresponding $\mbf Z/p^e\mbf Z$-torsors are split over the same vertices in this part of $\mbb T(\Omega)$. Thus, we will consider some currents which, over some "big" part of $\Omega$, coincides with the one corresponding to a homography invertible over $\omega$ and which is equivariant under some subgroup of finite index of $\Gal(\Omega/X)$. Then we consider the corresponding $\mbf Z/p^e\mbf Z$-torsor over some finite étale covering of $X$ which will behave like  $\mbf G_m\stackrel{z\mapsto z^{p^e}}{\to}\mbf G_m$ over some part of $\mbb T(\Omega)$.\\
We will deduce from this that the length of every loop of every finite topological covering of the graph of the stable model of $X$ can be recovered from the tempered fundamental group. A final combinatorial consideration will give us what we wanted.
\subsection{Preliminaries}
Let $K$ be a complete nonarchimedean field of characteristic 0, $K^\circ$
its integral ring, $k$ its residual field, assumed of characteristic
$p>0$, and let $\overline K$ be the completion of the algebraic closure of
$K$ (with integral ring $\overline K^{\circ}$).\\
From now on, we assume that we have chosen a compatible system of roots of 1 in $\overline K$, so that we may identify $\mu_n$ and $\mbf Z/n\mbf Z$. Thus we will often talk of $\mbf Z/n\mbf Z$-torsors over a curve over $\overline K$ when we should better talk of $\mu_n$-torsors.\\
Let $X_K$ a smooth curve of type $(g,n)$ over $K$ and let $Y_K\to X_K$ a
geometrically connected finite étale covering of $X_K$.\\

Let $\Pi=\gtemp(X_{\overline K})$ and $H\subset \Pi$ the open subgroup
of finite index corresponding to $Y_K$.\\

By changing $K$ by a finite extension, one may assume that
$X_K$ and $Y_K$ have a stable model over $K^{\circ}$, denoted $X$ and
$Y$.\\
One then have a unique morphism $\psi:Y\to X$ extending $Y_K\to X_K$.\\

By theorem~(\ref{mochi}), one may rebuild from the topological group $\Pi$ the semigraph
of groups $\mcal G$ of $X_{\overline k}$. Likewise, one may rebuild from $H$ the semigraph of groups $\mcal H$ of $Y_{\overline k}$.\\
We are now interested in reconstructing from $\gtemp(Y_K)\to\gtemp(X_K)$ what data can be recovered of the preimage of the cusps and vertices of the skeleton of $X_K$.\\

\subsubsection{Cusps}\label{cusp}
Let $y$ be a cusp of $Y_{\overline k}$ (corresponding to a
cusp $y'$ of $Y_{\overline K}$) and let $x$ be its image in
$X_{\overline k}$ (corresponding to a cusp $x'$ of $X_{\overline
    K}$).\\
Let $I\subset H^{(p')}$ be an inertia group of $y$. The image of $I$ in
$\Pi^{(p')}$ is an open subgroup (and thus it is nontrivial) of an
inertia group of $x$. Since the intersection of the inertia group of two different cusps is $\{1\}$, the image of $I$
is not contained in any other inertia group of a cusp of
$X_{\overline k}$, thus $x$ is characterized by the morphism $H\to \Pi$
as being the only cusp of $X_{\overline K}$ such that the
inertia groups of $y$ map by $H\to \Pi$ to inertia groups of $x$.\\
Thus one may rebuild from $H\subset\Pi$ the map from the
cusps of $Y_{\overline k}$ to the cusps of
$X_{\overline k}$.
In particular, one may decide whether the morphism $Y\to X$ is ramified at
$x'$.\\

\subsubsection{Vertices of the skeleton}\label{vertex}
If $X_0$ is an irreducible component of $X_{\overline k}$, let $Y_0$
be an irreducible component of $Y_{\overline k}$ which maps surjectively to
$X_0$. Then the morphism between the components of the
graphs of groups $\Pi^{(p')}_{Y_0}\to\Pi^{(p')}_{X_0}$ is open (in particular, its
image is non commutative) since it embeds in the following commutative diagram
\[\xymatrix{\Gal\big(\overline{K(Y_0)}/K(Y_0)\big) \ar@{->>}[d] \ar@{^{(}->}[r] &
  \Gal\big(\overline{K(X_0)}/K(Y_0)\big) \ar@{->>}[d] \\ \Pi^{(p')}_{Y_0} \ar[r] &
  \Pi^{(p')}_{X_0}},\]
where the upper arrow is an open embedding and the vertical arrows are projections.\\
Since $\Pi^{(p')}_{X_0}\to \Pi^{(p')}$ (defined up to conjugation) is injective, the image of
$\Pi^{(p')}_{Y_0}$ in $\Pi^{(p')}$ (defined up to conjugation) is non commutative, and thus
$\Pi_{X_0}^{(p')}$ is the only verticial subgroup of $\Pi^{(p')}$ which
contains the image of $\Pi^{(p')}_{Y_0}$.\\
Moreover, if $Y_0$ is an irreducible component of $Y_{\overline k}$ which
does not map surjectively onto an irreducible component of $X_{\overline k}$,
the image of $\Pi^{(p')}_{Y_0}$ in $\Pi^{(p')}$ is commutative, so the embedding $H\to\Pi$ decides which components of $Y_{\overline k}$ maps surjectively onto which components of
$X_{\overline k}$.\\
In particular, one may rebuild from $H\to \pi$ the cardinal of
the preimage by $Y_{\overline k}$ of the generic points of $X_{\overline k}$.\\
Equivalently, if $x_0$ is a vertex of the skeleton of
$X_{\overline K}^{\an}$, one may rebuild the cardinal of its preimage
in $Y_{\overline K}^{\an}$.\\
Indeed, let us write $\pi$ for the continuous map from the generic analytic fiber
to the special fibre. $x_1=\pi(x_0)$ is the generic point of an irreducible component of $X_{\overline k}$, every preimage $y_1$ of $x_1$
by $\psi$ is a generic point of an irreducible component of
$Y_{\overline k}$, $\pi^{-1}(y_1)$ is reduced to a single element
by~\cite[cor. 1.7]{berk2}, which must map to $x_0$ because
$\pi^{-1}(x_1)=\{x_0\}$ by~\cite[cor. 1.7]{berk2}. Thus
$\psi^{\an,-1}(x_0)$ is in natural bijection with $\psi^{-1}(x_1)$. 

\subsection{Case of $\mbf{P}^1\backslash\{z_1,\dots,z_n\}$}
Let $\sq=\alpha,\beta$.\\
Let $z_{\sq,1},\cdots,z_{\sq,n}\in\mbf Q_p^{\text{nr}}$, with
$n\geqslant 4$.\\
Let us write $X_{\sq}=\mbf{P}^1\backslash\{z_{\sq, 1},\dots,z_{\sq, n}\}$.\\
Let $\Pi_{\sq}=\gtemp(\mbf{P}^1\backslash\{z_{\sq, 1},\dots,z_{\sq,
  n}\})$. We already know that
an isomorphism $\phi:\Pi_{\alpha}\simeq\Pi_{\beta}$ induces an isomorphism
between the semigraphs of the stable reductions
of $\mbf{P}^1\backslash\{z_{\alpha,1},\dots,z_{\alpha,n}\}$. After
reordering the $z_{\beta, i}$, we may assume that this morphism of
semigraphs identifies the inertia subgroup (defined up to conjugation) of
the cusp $z_{\alpha,i}$ with the inertia subgroup of the cusp $z_{\beta,i}$.\\

\begin{thm}\label{thmp1} The isomorphism of graphs thus defined by $\phi$
  between the skeletons of
$(\mbf{P}^1\backslash\{z_{\alpha, 1},\dots,z_{\alpha, n}\})^{\an}$ and
$(\mbf{P}^1\backslash\{z_{\beta, 1},\dots,z_{\beta, n}\})^{\an}$ preserves
lengths of edges (\emph{i.e.} it induces an isomorphism of metric graphs). \end{thm}
Equivalently, for every $(i_1,i_2,i_3,i_4)$, one gets
the equality of cross-ratios' norms: 
\[|(z_{\alpha,i_1},z_{\alpha,i_2},z_{\alpha,i_3},z_{\alpha,i_4})|=|(z_{\beta,i_1},z_{\beta,i_2},z_{\beta,i_3},z_{\beta,i_4})|.\]
In fact, we will be able to prove this result without assuming $z_{\sq,1},\cdots,z_{\sq,n}\in\mbf Q_p^{\text{nr}}$ after studying the case of an elliptic curve (this will give the result for $p\neq 2$) and without any assumption after studying the case of a Mumford curve.\\

\dem According to section~\ref{cusp}, if $H_{\alpha}$ and $H_{\beta}=\phi(H_{\alpha})$
are finite index open subgroups corresponding to each other by $\phi$, we know
that the corresponding covering of $X_\alpha$ is ramified outside
$\{z_{\alpha,i_1},z_{\alpha,i_2},z_{\alpha,i_3},z_{\alpha,i_4}\}$ if and
only if the corresponding covering of $X_\beta$ is ramified outside
$\{z_{\beta,i_1},z_{\beta,i_2},z_{\beta,i_3},z_{\beta,i_4}\}$.\\
Thus $\phi$ induces an isomorphism between the
$\gtemp(\mbf{P}^1\backslash\{z_{\sq, i_1}, z_{\sq, i_2}, z_{\sq, i_3},
z_{\sq, i_4}\})$ (it is the projective limit of the discrete quotients
$G_{\sq}$ of $\Pi_{\sq}$ which have a finite quotient $G_{1,\sq}$
corresponding to a covering of $\mbf P^1$ unramified outside
$\{z_{\sq,i_1},z_{\sq,i_2},z_{\sq,i_3},z_{\sq,i_4}\}$ such that
$\Ker(G_{\sq}\to G_{1,\sq})$ is a free group (or, equivalently here, is
torsion free)).\\

We thus may assume that $n=4$. We may also assume, after a base change of
an homography, that the cusps $0,1,\infty$ are
$\lambda_{\sq}$. Moreover, we may assume that $X_{\sq}$ does not have good
reduction (in which case there is nothing to prove), and thus that $v_p(\lambda_{\sq}-1)>0$
(it is an integer since $\lambda\in\mbf Q_p^{\text{nr}}$) after permuting
$0,1$ and $\infty$ by another homography. We now have to prove that
$v_p(\lambda_{\alpha}-1)=v_p(\lambda_{\beta}-1)$.\\
Assume, ab absurdo, that
$v_p(\lambda_{\alpha}-1)<v_p(\lambda_{\beta}-1)$.
Let $e:=v_p(\lambda_{\beta}-1)-1\geqslant v_p(\lambda_{alpha}-1)$ and let
$H_{\sq}$ be the subgroup of
$\Pi_{\sq}$ of index $p^e$
corresponding to the unique connected covering $Y_{\sq}\to X_{\sq}$
of degree $p^e$ unramified outside $0$ and $\infty$ (this is the morphism $z\mapsto z^{p^e}$
from $\mbf P^1$ to $\mbf P^1$). One has $\phi(H_\alpha)=H_\beta$ by section~\ref{cusp}.\\
But, according to lemma~(\ref{decompositionpuissance}), $B(1,r)$ (the point of the Berkovich space of $\mbf P^1$ corresponding to
the ball of radius $r$ and center $1$) has $p^{e-1}$ preimages if
$p^{-(e+\frac{1}{p-1})}\leqslant r <p^{-(e-1+\frac{1}{p-1})}$ and has $p^{e}$
preimages if $r<p^{-(e+\frac{1}{p-1})}$.\\
Thus if $p\neq 2$, $B(1,|\lambda_{\alpha}-1|)$ has $p^{e-1}$ preimages if
$Y_\alpha$ and $B(1,|\lambda_{\beta}-1|)$ has $p^e$ preimages in $Y_\beta$,
which contradicts the result of part~\ref{vertex}.\\

In the case $p=2$ and $e\geqslant 2$,  $B(1,|\lambda_{\alpha}-1|)$ has $2^{e-2}$ preimages in
$Y_\alpha$ and $B(1,|\lambda_{\beta}-1|)$ has $2^{e-1}$ preimages in $Y_\beta$,
which contradicts the result of part~\ref{vertex}.\\
If $e=1$, and therefore $v_2(\lambda_{\alpha})=1$ end $v_2(\lambda_{\beta})=2$, the
semigraph of the reduction of $Y_{\alpha}\to X_{\alpha}$ is (we marked the
different cusps of $Y_{\alpha}$, and
$\sqrt{\lambda_{\sq}}$ is a square root of $\lambda_{\sq}$):  
\[\xymatrix{ \infty & & & \sqrt{\lambda_{\alpha}} \\  & & \sommet \ar@{-}[ur] \ar@{-}[r] & -\sqrt{\lambda_{\alpha}}\\ \sommet
  \ar@{-}[uu] \ar@{-}[dd] \ar@{-}[r] & \sommet \ar@{-}[ur] \ar@{-}[dr]&  & \\ & & \sommet
  \ar@{-}[r] \ar@{-}[dr] & 1\\ 0 & & & -1}\]
The semigraph of the stable reduction of $Y_{\beta}$ is 
\[\xymatrix{ \infty & & & \sqrt{\lambda_{\beta}} \\ & & & -\sqrt{\lambda_{\beta}}\\ \sommet \ar@{-}[uu] \ar@{-}[dd] \ar@{-}[rr] & & \sommet \ar@{-}[uur]
  \ar@{-}[ur] \ar@{-}[dr] \ar@{-}[ddr] & \\ & & & 1\\  0 & & & -1}\]
They are not isomorphic, and thus one also gets a contradiction.\findem

\subsection{Case of a punctured elliptic curve}
Let $\sq=\alpha,\beta$\\
Let $X_{\sq}$ be two punctured Tate curves $\mbf C_p^*/q_{\sq}^{\mbf
  Z}-\{1\}$ with $|q_{\sq}|<1$.\\
Let $\Pi_{\sq}=\gtemp{X_{\sq}})$ and let $\phi:\Pi_\alpha\simeq\Pi_\beta$
be an isomorphism. 
\begin{thm}
$|q_{\alpha}|=|q_{\beta}|$\label{ellcurve}
\end{thm}
\dem\\
Let us choose integers $n$, $l$ and $m$ such that:
\begin{itemize}
\item $n$ is prime to $p$ and$n\geqslant \frac{v_p(q_\beta)v_p(q_\alpha)(p-1)}{|v_p(q_\beta)-v_p(q_\alpha)|p}$,
\item $l\geqslant 1+\frac{2np}{(p-1)v_p(q_{\sq})}$,
\item $m\geqslant \frac{2l}{n}$.
\end{itemize}
Let  $H_{0,\sq}=[\Pi_{\sq},\Pi_{\sq}]\Pi_{\sq}^n$ be the preimage in
$\Pi_{\sq}$ of the image by the multiplication by $n$ in the abelianized
group of $\Pi_{\sq}$. $\phi$ induces an isomorphism
$H_{0,\alpha}\to H_{0,\beta}$.\\
The covering $Y_{0,\sq}$ of $X_{\sq}$ corresponding to $H_{\sq}$ is the multiplication
$\overline X_{\sq}\stackrel{\times n}{\to}\overline X_{\sq}$
by $n$ on the elliptic curve $\overline{X}_{\sq}$.\\

Let now $H_{1,\sq}$ be the subgroup of $H_{0,\sq}$ corresponding to
the unique connected topological covering $Y_{1,\sq}$ of degree $n$ of
$Y_{0,\sq}$.\\
$Y_{1,sq}\simeq \mbf C_p^*/q^{m\mbf
  Z}-\{q^{\frac{a}{n}}\zeta^b\}_{(a,b)\in\mbf Z^2}$ où $q^{\frac{1}{n}}$
is an $n$-th root of $q$ and $\zeta$ is an $n$-th root of $1$.\\ 
The semigraph of the stable reduction of $Y_{1,\sq}$ has $mn$ vertices
joined in circle (the distance of such two successive vertices is
$v_p(q_{\sq})/n$) and  $n$ cusps abuts to each vertex.\\
$\phi$ induces an isomorphism $H_{1,\alpha}\simeq H_{1,\beta}$ which itself
induces an isomorphism between the semigraphs of the stable reduction of
$Y_{1,\sq}$. Let us number from $0$ to $mn-1$ the vertices of the graphs by
following the circle compatibly with the isomorphism induced by $\phi$ (let us write $x_{\sq,0},\cdots, x_{\sq,mn-1}$ for the
corresponding vertices of the skeleton of $Y_{1,\sq}$).\\
Let $z_{\alpha,0}$ and $z_{\alpha,,l}$ be two cusps of
$Y_{1,\alpha}$ abutting to the vertices of the graph numbered $0$ and $l$
respectively. Let $z_{\beta,0}$ and $z_{\beta,1}$ be the corresponding
cusps of $Y_{1,\beta}$.\\

Let us now focus on $\mbf F_p$-torsors over $\overline Y_{1,\sq}$ which are
unramified outside $z_{\sq,0}$ and $z_{\sq,l}$. They are the elements of an
$\mbf F_p$-vector space $V_{\sq}$ of dimension $3$.\\
Recall that we chose $m,n$ and $l$ so that $\frac{l-1}{n}v_p(q_{\sq})>\frac{2p}{p-1}$ and
$\frac{mn-l}{n}v_p(q_{\sq})>\frac{l}{n}v_p(q_{\sq})$.\\
We will now describe a base of this $\mbf F_p$-vector space.\\

Let $S_{\sq}$ be the universal topological covering of $\overline
Y_{1,\sq}$, we identify it with $\mbf P^1\backslash\{ 0,\infty
\}\subset\mbf P^1$. Let
$s_{\sq,0}$ and $s_{\sq,l}$ the unique preimages in $S_{\sq}$ of $Z_{\sq,0}$ and
$z_{\sq,l}$ of norm $1$ and $|q|^{\frac{l}{mn}}$ and let $U_{1,\sq}\subset
S_{\sq}$ be the open crown $\{|q_{\sq}|^{\frac{l-mn}{2n}}>|z|>|q_{\sq}|^{\frac{l+mn}{2n}}\}$
and $U_{2,\sq}\subset S_{\sq}$ the open annulus
$\{|q_{\sq}|^{\frac{l}{n}}p^{-\frac{p}{p-1}}>|z|>|q_{\sq}|^mp^{\frac{p}{p-1}}\}$. The
maps fom $U_{1,\sq}$ and from $U_{2,\sq}$ to $\overline Y_{1,\sq}$  still
are open embeddings which, together, cover $\overline Y_{1,\sq}$.\\
Let $T_{1,\sq}$ the restriction to $U_{1,\sq}$ of the ramified
(only over $s_{\sq,0}$ and $s_{\sq,l}$) covering $\mbf
P^1\to\mbf P^1:z\mapsto \frac{s_{\sq,0}z^p+s_{\sq,l}}{z^p+1}$, which is Galois of
Galois group isomorphic to $\mbf Z/p\mbf Z$, and let us choose such an
isomorphism to get a $\mbf Z/p\mbf Z$-torsor. Let $T_{2,\sq}$ be the
trivial $\mbf Z/p\mbf Z$-torsor over $U_{2,\sq}$  and let $T_{3,\sq}=T_{1,\sq}\coprod
T_{2,\sq}\to U_{3,\sq}=U_{1,\sq}\coprod U_{2,\sq}$.\\
Over $U_{1,\sq}\times_{\overline Y_{1,\sq}} U_{2,\sq}$ (which has two
connected components), $T_{1,\sq}$ is trivial. Let us choose a 
trivialization : one may now descend $T_{3,\sq}\to U_{3,\sq}$ into
a $\mbf Z/p\mbf Z$-torsor $T_{\sq}\to \overline Y_{1,\sq}$, which is only ramified
over $z_{\sq,0}$ and $z_{\sq,l}$.\\
According to lemma~(\ref{decompositionpuissance}), $T_{\sq}\to \overline Y_{1,\sq}$ is split over $x_i$ (with $i\in [0,mn-1]$)
if and only if
$|q_{\sq}|^{\frac{i}{n}}\in
[|q_{\sq}|^mp^{\frac{p}{p-1}},|q_{\sq}|^{\frac{l}{n}}p^{-\frac{p}{p-1}}]$,
\emph{i. e.} \[i\in I_{1,\sq}:=[l+\frac{np}{v_p(q_{\sq})(p-1)},mn-\frac{np}{v_p(q_{\sq})(p-1)}].\]
There is such an integer thanks to the assumption about $l,m$ and $n$ because
\[\lg(I_{1,\sq})=mn-l-2\frac{np}{v_p(q_{\sq})(p-1)}\geqslant 1\] (where $\lg$
denotes the length of an interval).\\

Likewise, let $s'_{\sq,l}$ be a preimage of $z_{\sq,l}$ of norm
$|q|^{\frac{l-mn}{n}}$, let $U'_{1,\sq}$ be the crown
$\{|q_{\sq}|^{\frac{l-2mn}{2n}}>|z|>|q_{\sq}|^{\frac{l}{2n}}\}$,
$U'_{2,\sq}$ the crown $\{p^{-\frac{p}{p-1}}>|z|>|q_{\sq}|^{\frac{mn-l}{n}}p^{\frac{p}{p-1}}\}$.
They are open subsets of $\overline Y_{1,\sq}$ and cover it.\\
Let $T'_{1,\sq}$ be a $\mbf Z/p\mbf Z$-torsor over $U'_1$ obtained as the 
restriction of a $\mbf Z/p\mbf Z$-torsor over $\mbf P^1$ ramified only above
$z_{\sq,0}$ and $z'_{\sq,l}$. It is trivial over $U'_{1,\sq}\cap
U'_{2,\sq}$, and, by choosing such a trivialization, one gets by descent a
$\mbf Z/p\mbf Z$-torsor $T'_{\sq}$ over $\overline Y_{1,\sq}$,
ramified only over $z_{\sq,0}$ and $z_{\sq,l}$.\\
$T_{\sq}\to \overline Y_{1,\sq}$ is split over $x_i$ (with $i\in [0,mn-1]$)
if and only if
$|q_{\sq}|^{\frac{i}{n}}\in
[|q_{\sq}|^{\frac{l}{n}}p^{\frac{p}{p-1}},p^{-\frac{p}{p-1}}]$, \emph{i. e.}
\[i\in
I_{2,\sq}:=[\frac{np}{v_p(q_{\sq})(p-1)},l-\frac{np}{v_p(q_{\sq})(p-1)}].\]
There is such an integer thanks to the assumption about $l,m$ and $n$
because \[\lg(I_{2,\sq})=l-2\frac{np}{v_p(q_{\sq})(p-1)}\geqslant 1\].\\

Let finally $T''_{\sq}$ be the essentially unique connected topological covering
of degree $p$ of $\overline Y_{\sq}$ and let us choose an isomorphism from $\mbf
Z/p\mbf Z$ onto its Galois group so that it becomes a $\mbf Z/p\mbf
Z$-torsor. Let us then show that $T_{\sq},T'_{\sq},$ and $T''_{\sq}$
constitute a base of $V_{\sq}$.\\
Let $i$ be an integer in $I_{1,\sq}$. As $I_{1,\sq}\cap
I_{2,\sq}=\emptyset$ and as $T''$
is everywhere split, if $aT_{\sq}+bT'_{\sq}+cT''_{\sq}$ (the linear
combination is in the sense of the structure of vector space $V_{\sq}$) is
split over $x_i$, $c=0$. By the same argument with $I_{2,\sq}$, if
$aT_{\sq}+bT'_{\sq}+cT''_{\sq}=0$, one gets $b=0$, but as $T_{\sq}$
is not trivial either, one indeed gets that $T_{\sq},T'_{\sq}$ and
$T''_{\sq}$ constitute a base of $V_{\sq}$.\\

Assume now ab absurdo that $|q_\alpha|>|q_\beta|$, then $I_{1,\alpha}\subset
I_{1,\beta}$.\\
Let $T_{0,\beta}=aT_{\beta}+bT'_{\beta}+cT''_\beta$ be the image of $T_{\alpha}$ by $(\phi^{-1})^*:V_{\alpha}\to
V_{\beta}$ (indeed, $phi$ induces such a $(\phi^{-1})^*$ according to part~\ref{cusp}).\\
Let $i\in I_{1,\alpha}$, $T_{\alpha}$ is split over $x_{\alpha,i}$,
so, according to section~\ref{vertex}, $T_{0,\beta}$ is split over $x_{\beta,i}$,
and thus $c=0$. Thus, $T_{0,\beta}$ is split over every $x_i$ if
$b=0$ or exactly over the $x_i$ with $i\in I_{1,\beta}$ otherwise.\\
Yet, according to section~\ref{vertex}, $T_{0,\beta}$ must be split exactly over
the $x_{\beta,i}$ such that $T_{\alpha}$ is split over $x_{\alpha,i}$,
\emph{i. e.} $i\in I_{1,\alpha}$.\\
Thus, one cannot be in the case $b=0$ and  $I_{1,\beta}$ and $I_{1,\alpha}$
must contain exactly the same integers.\\
But this cannot be if  $n\geqslant
\frac{v_p(q_\beta)v_p(q_{\alpha})(p-1)}{(v_p(q_\beta)-v_p(q_\alpha))p}$,
because
\[\lg(I_{1,\beta})-\lg(I_{1,\alpha})=2\frac{np}{v_p(q_{\alpha})(p-1)}-2\frac{np}{v_p(q_{\beta})(p-1)}\geqslant
2.\]   
\findem

\emph{Remark:}\\
Assume $p\neq 2$. Let $\{z_{\sq,1},\cdots, z_{\sq,4}\}$ be
  four rational points in
  $\mbf P^1_{\mbf C_p}$, and let $\phi$ be an isomorphism between the
  $\gtemp(\mbf P^1\backslash\{z_{\sq,1},\cdots, z_{\sq,4}\})$, which
  identifies $z_{\alpha,i}$ to $z_{\beta,i}$ by the identification of the
  cusps of the graphs of the stable reduction (we assume that the curves
  have bad reduction, and that the length of the edge of the graph is
  $l_{\sq}$). Let $E_{\sq}$ be the unique
  $\mbf Z/2\mbf Z$ covering of $\mbf P^1$ ramified over
  $\{z_{\sq,1},\cdots, z_{\sq,4}\}$ and only over those points (the
  subgroups of index 2 corresponding to $E_{\sq}$ thus maps to each other). $E_{\sq}$
  is a Tate curve with $|q_{\sq}|=2l_{\sq}$. By theorem~(\ref{ellcurve}),
  $|q_{\alpha}|=|q_{\beta}|$ so $l_{\alpha}=l_{\beta}$, which proves
  again theorem~(\ref{thmp1}) without assuming that the points are in
  $\mbf Q_p^{\text{nr}}$.\\
If $p=2$ and $l_{\sq}>4$, then $E_{\sq}$ is also a Tate curve, $|q_{\sq}|=2l_{\sq}-8$ and
the previous argument works as well.

\subsection{Case of a Mumford curve}
\subsubsection{Reminder on Mumford curves and currents}
Let $X$ be a Mumford curve of genus $g$ over $\overline K$, let  $\Omega\subset \mbf P^1$
be its universal topological covering, and $\Gamma=\Gal(\Omega/X)$, so that $X=\Omega/\Gamma$.\\
Let $\Phi$ be the retraction of $\Omega$, as a Berkovich space, onto
its skeleton $\mbb T=\mbb T(\Omega)$.\\
For $z,z'\in \Omega$, let us set\[d(z,z')=\sup_{x_1,x_2\in\mbf
  P^1\backslash\Omega}|v_p(\frac{z'-x_1}{z-x_1}\frac{z-x_2}{z'-x_2})|\]
which is invariant under homographies. Moreover $d(z,z')$ only depends on
$\Phi(z)$ and $\Phi(z')$, and it is nothing else than the distance between $\Phi(z)$
and $\Phi(z')$ for the usual metric structure on the tree of $\Omega$ (that we will also
denote by $d$, so that $d(z,z')=d(\Phi(z),\Phi(z'))$).\\
If $z\in\Omega$ and $\lambda>0$, let $U_{z,\lambda}$ denote $\{z\in\Omega/d(z,z')\leqslant \alpha$\}. It is an affinoid subspace of $\Omega$.\\

Let $\mcal L=\mbf P^1\backslash \Omega$, a compact subset of $\mbf P^1$.
In~\cite[1.8.9]{FvdP2}, Fresnel and van der Put define a \emph{measure} on a profinite topological space $Z$ as a function $\mu:\{\text{compact open subsets of Z}\}\to\mbf Z$ such that $\mu(U_1\cup U_2)+\mu(U_2\cap U_1)=\mu(U_1)+\mu(U_2)$ for every compact open subsets $U_1,U_2$ of $Z$. The group of measures on $Z$ such that $\mu(Z)=0$ is then denoted by $M_0(Z)$.\\
One can then associate to $f\in
O(\Omega)^*$ a measure $\mu_f$ on $\mcal L$. Following~\cite{FvdP2}, we will also call hole of an affinoid $U$ subspace of $\mbf P^1$ every connected component of the complement of this affinoid, and will note $t(U)$ the set of its holes.\\
One gets the following exact sequence of groups~(\cite[Prop. 1.8.9]{FvdP2}):
\[1\to \overline K^*\to O(\Omega)^*\to M_0(\mcal L)\to 0\]
More precisely, according to~\cite[1.8.10 ex. $\beta$]{FvdP2}, if $a,b\in\mcal L$ and $f:z\mapsto\frac{z-a}{z-b}$, $\mu_f=\delta_a-\delta_b$.\\
For a general $f_0\in O(\Omega)^*$ and $\mu=\mu_{f}$, $\mu$ is a weak limit of a sequence  of $\mbf Z$-linear combinations $(\mu_k)_{k\in\mbf N}$ of Dirac measures. Then, if $\mu_k=\sum n_i(k)\delta_{a_i(k)}$, let $f_k=\prod (1-\frac{a_i(k)}{z})^{n_i(k)}$. Then $\mu_{f_k}=\mu_k$, and $(f_k)$ tends uniformly on every affinoid to an invertible function $g$. Then $\mu_g=\mu$ and thus $g=\lambda f$ with $\lambda \in \overline K^*$.\\

According to~\cite[prop. 1.1]{vdp1}, one also has the following exact sequence (where $C(\mbb T)$ denotes the group of currents with integral coefficients on $\mbb T$):
\[1\to \overline K^*\to O(\Omega)^*\to C(\mbb T)\to 0\]
which gives, with the previous one an isomorphism $M_0(\mcal L)\to C(\mbb T)$.\\
One can describe this isomorphism in the following way : if $e$ is an oriented edge of $\mbb T$, $\mbf P^1\backslash\Phi^{-1}(e)$ has two connected components, and $\Phi^{-1}(e)\cap\mcal L=\emptyset$ so one gets a partition of $\mcal L$ by two open subsets $\mcal L_1(e)$ at the beginning of $e$ and $\mcal L_2(e)$ at the end of $e$ (one has for $\mu\in M_0(\mcal L)$, $\mu(\mcal L_1(e))=-\mu(\mcal L_2(e))$. Then $\mcal C(e)=\mu(\mcal L_2(e))$.\\
more generally, $\mbb K$ is a finite connected subgraph of $\mbb T$ (containing at least one edge), $\Phi^{-1}(\mbb K)$ is an affinoid contained in $\Omega$. There is a natural bijection between $t(\Phi^{-1}(\mbb K))$ and the set of end points of $\mbb K$. Then, if $\mu\in M_0(\mcal L)$, $C\in t(\Phi^{-1}(\mbb K))$ and $e$ is the only oriented edge of $K$ ending at the end point of $\mbb K$ corresponding to $C$, then $C\cap \mcal L=\mcal L_2(e)$ and thus $\mcal C(e)=\mu(C\cap \mcal L)$. In particular $\mu(C\cap \mcal L)=0$ for every $C\in t(\Phi^{-1}(\mbb K))$ if and only if $\mcal C(e)=0$ for every terminal edge of $\mbb K$, if and only if $\mcal C$ is zero on $\mbb K$.\\

Let $\Theta$ denote the group of theta functions of $X$, that is the group of $f\in O(\Omega)^*$ such that for every $\gamma\in \Gamma$, $z\mapsto f(\gamma z)/f(z)$ is a constant function (this means that the corresponding current is $\Gamma$-equivariant). Then, one has the following exact sequence :
\[1\to \overline K^*\to \Theta \to C(\Gamma)\to 0\]
and thus $\Theta/\overline K^*$ is a free $\mbf Z$-module of rank $g$.\\
One deduces from~\cite[prop. 2.1]{vdp2} that :
\begin{prop} for every $n\geq 2$ and every $\mbf Z/n\mbf Z$-torsor $Y$ over $X$,
there exists an element $\theta$ in $\Theta$, unique modulo $\overline K^*\Theta^n$ such that $Y\times_X \Omega=\Omega[f]/(f^n=\theta)$ where $\Omega[f]/(f^n=\theta)$ denotes the pullback of the $\mbf Z/n\mbf Z$-torsor $\mbf G_m \stackrel{z\mapsto z^n}{\to}\mbf G_m$ along $\Omega\stackrel{\theta}{\to}\mbf G_m$.\\
Conversely, for every $\theta$ in $\Theta$, there exists a $\mbf Z/n\mbf Z$-torsor $Y$ over $X$ such that $Y\times\Omega=\Omega[f]/(f^n=\theta)$.\end{prop}
\dem Follow the notations of~\cite[section 2]{vdp2}. Then $\Omega_*$ is a connected component of $Y\times_X\Omega$.\\
Suppose first that $Y\times_X\Omega$ is connected, so that $\Omega_*=Y\times_X\Omega$. Then, according to~\cite[prop. 2.1]{vdp2} there is a unique lattice $T$ in $(\Theta/\overline K^*)\otimes \mbf Q$ containing $\Theta/\overline K^*$ such that, as a $X$-covering, $\Omega_*=\Omega(T):=\Omega\times_{\Spec \overline K[\Theta/\overline K^*]}\Spec \overline K[T]$. Then $T/(\Theta/\overline K^*)$ is isomorphic to $\mbf Z/n\mbf Z$, and choosing a $\bar f$ generator of $T/(\Theta/\overline K^*)$) amounts to choosing a $\mbf Z/n\mbf Z$-torsor structure on $\Omega_*$.\\
If one takes $\bar f^i$ with $i\in (\mbf Z/n\mbf Z)^*$ as another generator, the torsor one gets is $\Omega_*^{\otimes i}$. Thus by changing the generator, one can get all the $\phi(n)$ $\mbf Z/n\mbf Z$-torsor structure on the covering $\Omega_*$ of $\Omega$.\\
If $\bar f$ is the generator corresponding to the $\mbf Z/n\mbf Z$-torsor structure on $Y\times_X\Omega$, then $\theta$ is any lifting of $\bar f^n\in T^n/(\Theta/\overline K^*)^n$.\\
In the general case $\Omega_*$ acquires the structure of a $\mbf Z/m\mbf Z$-torsor over $\Omega$ with $m|n$, and as before one can find a unique $\theta_0$ modulo $\overline K^*\Theta^m$ such that $\Omega_*=\Omega[f]/(f^m=\theta_0)$. Then $Y\times_X\Omega=\Ind_{\mbf Z/m\mbf Z}^{\mbf Z/n\mbf Z}\Omega_*$, and thus $Y\times_X\Omega=\Omega[f]/(f^n=\theta_0^{m/n})$, and $\theta=\theta_0^{m/n}$.\\
The second statement comes from the fact that, if $\Omega_*$ is a connected component of $\Omega[f]/(f^n=\theta)$, $\Gal(\Omega_*/X)$ is (non canonically) isomorphic to the direct product of $\Gal(\Omega_*/\Omega)$ and $\Gamma$ (according to~\cite[section 2, intro.]{vdp2}). Thus $\Omega_*$ can descend (non canonically) to $X$ by considering $Y_0=\Omega_*/N$ where $N$ is some complement of $\Gal(\Omega_*/\Omega)$ and $\Gal(\Omega_*/X)$ (and thus $\Omega[f]/(f^n=\theta)$ by taking a direct sum of $Y_0$s.\findem

\emph{Remark:} One could also show it by considering $\widetilde J=\mcal Hom(\Theta/\overline K^*,\mbf G_m)\to J$ where $J$ is the Jacobian variety of $X$ and $\widetilde J$ is its universal topological covering, and by showing that $\ga(\widetilde J)$ is a direct summand of $\ga(J)$.

\subsubsection{Preliminary results on the ramifications of torsors corresponding to currents}
\begin{prop} Let $z\in\Omega,\lambda>0$. Let $f\in O(\Omega)^*$ such that $f(z)=1$.
  Let $\mu$ be the measure on $\mcal L$ corresponding to $f$
  and let us assume $\mu(C\cap \mcal L)=0$ for every hole of $U_{z,\lambda}$.\\
Then $\forall z'\in U_{z,\lambda}$, $|f(z')-1|\leqslant
p^{d(z,z')-\lambda}$\end{prop}
\dem To simplify, assume $z=\infty$.\\
According to~\cite[1.8.10 ex. $\beta$]{FvdP2} and~\cite[prop. 1.8.9.(i)]{FvdP2},
$f=\lim f_k$ uniformly on every affinoid of $\Omega$ (in particular over $U_{z,\lambda}$) where $f_k$ is of the following form
\[f_k(z')=\prod^{s_k}_{i=1} (1-\frac{x_{i,k}}{z'})^{n_{i,k}},\]
and $\mu_k=\sum_i n_{i,k}\delta_{a_i(k)}\to \mu$.\\
For $k$ big enough, $\mu_k(C)=0$ for every hole $C$ of
$U_{z,\lambda}$ and $|f-f_k|_{U_{z,\lambda}}\leqslant p-{\lambda}$.\\
We thus only have to prove the result for $f_k$, which is a product of function such as $g:z'\mapsto
z'=(1-\frac{x_1}{z'})(1-\frac{x_2}{z'})^{-1}=\frac{z'-x_1}{z'-x_2}$ with $x_1$
and $x_2$ in the same hole $C$ of $U_{z,\alpha}$.\\
We thus only have to prove the result for $g$, which is easily seen.\findem
\begin{cor} Let $f\in O(\Omega)^*$ such that $f(z)=1$. Let $U$ an affinoid
  of $\Omega$ such that for every hole $C$ of $U$, $\mu(C)=0$.\\
Assume that, for all $z''\in\Phi^{-1}([\Phi(z),\Phi(z')]),
U_{z'',\lambda}\subset U$.\\
Then $|f(z')-1|\leqslant p^{-\lambda}$.\end{cor}
\dem Let $\epsilon>0$. Let $(z_i)_{i=0\cdots n}$ such that $z_0=z$, $z_n=z'$, $\Phi(z_i)\in
[\Phi(z),\Phi(z')]$ and $d(\Phi(z_i),\Phi(z_{i+1}))\leqslant \epsilon$.\\
Then, according to the previous proposition
$|\frac{f(z_{i+1})}{f(z_i)}-1|\leqslant p^{\epsilon-\lambda}$.\\
Thus $|f(z'')-1|\leqslant \sup|f(z_{i+1})-f(z_i)|\leqslant
p^{\epsilon-\alpha}$. Ones get the result by letting $\epsilon$ going to 0.\findem
\begin{cor}\label{totdecrev} Let $f$ as previously, $U$ such that $\mu(C)=0$ for every hole $C$ of $U$. Let $e$ a positive integer. Let $\lambda>e+\frac{1}{p-1}$, let $Y\to
  \Omega$ be finite covering obtained by pulling back $\mbf G_m\stackrel{z\mapsto z^{p^e}}{\to}\mbf G_m$ along $f:\Omega\to \mbf G_m$. Let
  $V\subset U$ be such that $\forall z\in V, U_{z,\lambda}\subset U$.\\
Then $Y$ is split over $V$.\end{cor}
\dem We may assume $V$ connected, because we only have to prove the result for every connected component of $V$. Let $z\in V$. Multiplying $f$ by a constant, which do not change $Y$, we may assume that $f(z)=1$.\\
From the previous corollary, $f(V)\subset
D(1,p^{-\lambda})$. But, according to lemma~(\ref{decompositionpuissance}), $\mbf G_m\stackrel{z\mapsto z^{p^e}}{\to}\mbf G_m$
is split over $D(1,p^{-\lambda})$, which ends the proof.\findem

\begin{prop}\label{nondecrev} Let $\mcal C$ a current on $\mbb T(\Omega)$, corresponding
  to an invertible function $f$ on $\Omega$. Let $a$ a vertex of $\mbb
  T$ such that the restriction $\mcal C_a$ of $\mcal C$ the star of $a$ is not zero modulo $n$. Then, if $Y\to
  \Omega$ is the finite covering obtained by pullback of $\mbf G_m\stackrel{z\mapsto z^n}{\to}\mbf G_m$ along $f:\Omega\to  \mbf G_m$ , $Y\to
  \Omega$ is not split over $a$.\end{prop}
\dem $Y\to \Omega$ is split over $a$ if and only if there exists $f_1\in \mcal O_{\Omega,a}$ such that $f_1^n=f_{|\mcal  O_{\Omega,a}}$.\\
By multiplying $f$ by a constant, we may assume $|f|_a=1$.\\
If $f_1$ exists, by looking at the residue field $\tilde
k_{\Omega,a}\simeq\tilde k(X)$, $\tilde f_1^n=\tilde f$. If $\tilde
f_1(z)=\lambda\prod(z-a_i)$, $\tilde f=\lambda^n\prod(z-a_i)^n$, and so all the poles and zeros are of order a muliple of $n$, which ends the proof.\findem

\subsubsection{The metric graph of the reduction and the tempered fundamental group}
Assume now that $X_\alpha$ and $X_\beta$ soient are two
Mumford curves over $\overline K$, but are pullback of curves over $K$ (so that we may use~\cite[ex. 3.10]{mochi}), and that there is an isomorphism $\phi:\gtemp(X_\alpha)=\gtemp(X_\beta)$,
which thus inducts an isomorphism of graphs $\mbb G(X_\alpha)\to \mbb
G(X_\beta)$. Hence an isomorphism $\mbb T(\Omega_{\alpha})\simeq
T(\Omega_{\beta})$.\\
\begin{thm}\label{mumford}The isomorphism $\mbb G_{\alpha}\to\mbb G_{\beta}$ of graph induced by an isomorphism $\gtemp(X_\alpha)\to\gtemp(X_\beta)$ of topological groups is in fact an isomorphism of metric graphs.\end{thm}

\emph{Remark:} Suppose $X_\alpha$ and $X_\beta$ are projective lines minus four points $a_{\sq},b_{\sq},c_{\sq},d_{\sq}$ with $v_p(a_{\sq},b_{\sq},c_{\sq},d_{\sq})=e_{\sq}$with an isomorphism $\phi$ between their tempered fundamental groups, and suppose $p\neq 3$. One can consider a covering $X'_\alpha$ of order three of $X_\alpha$ such that the restriction to each irreducible component of the stable reduction is connected but it is split over the double point of the stable reduction (and let $X'_\beta$ be the covering of $X_\beta$ corresponding to $X'_\alpha$ by $\phi$, it satisfies the same properties). Then $\overline X'_{\sq}$ is a Mumford curve (the corresponding covering of each irreducible component of the stable reduction of $X_{\sq}$ is ramified in only two points, so it is a covering by a projective line) whose tree has three edges, each of length $e_{\sq}$. Thus, by what has just been proven  $e_\alpha=e_\beta$, which ends the proof for punctured lines.\\
One could in fact do quite the same for more general punctured Mumford curve, by considering, for an edge $e$ of the graph, a tamely ramified covering by a Mumford curve such that there is an edge $e'$ over $e$ which is actually an edge of the graph of the stable reduction of the compactification of this covering.\\

\dem 
A \emph{loop} of a graph is a cyclic sequence of oriented edges of the graph such that the end of an edge is the beginning of the following edge, and which never goes by the same vertex or (unoriented) edge twice.\\
Let $C$ be a loop of this graph, and note $\lg_{\sq}(C)$ the length
of $C$ with respect to the metric $d_{\sq}$ on $\mbb G$. Let $\widetilde C$ be the universal covering of $C$, let $\widetilde C\to\mbb T$ a lifting
of $\widetilde C\to \mbb G$  and let
$z_0$ a vertex of $\mbb T$ which belongs to
$\widetilde C$, let us number then $(z_i)_{i\in \mbf Z}$ the vertices of
$\widetilde C$ with the same image that $z_0$ in $\mbb G$. Let $L$ be another loop (we choose an orientation on it) of $\mbb G$ (there must be another loop since $g>1$), let $\widetilde L$ be a lifting of $L$ to
$\mbb T$, let $r_{\sq}=d_{\sq}(\widetilde L,z_0)$ (we may assume, by changing the numbering of the $z_i$ that
for $n\geqslant 0$, $d_{\sq}(\widetilde L,z_n)=r_{\sq}+n\lg_{\sq}(C)$) and
let $z'_{0}$ the point of $\widetilde L$ the nearest from $z_0$ (this does not depend on $\sq$).\\

Let us note, for $z$ a vertex of $\mbb T$, $F_z$ the connected component of
$\Omega\backslash\{\text{open edges of } \widetilde L\}$ which contains $z$.\\

Let $e\geqslant 1$ be an integer.\\
Let $K_0$  be a connected finite subgraph of $\mbb T$ containing $z_1$ (so that $\Phi^{-1}(|K_0|)$ is compact by properness of $\Phi$).\\
Let $K_{\sq}$ a connected finite subgraph (coinciding for $\sq=\alpha$ and $\sq=\beta$) $\mbb T(\Omega_{\sq})$ such that for every $z\in \Phi^{-1}(|K_0|)$
$U_{z_,e+2}\subset \Phi^{-1}(|K|)$ and such that, for every vertex $z$ of $\widetilde L\cap K_{\sq}$, $\{z'\in
F_z|d(z',z)\leqslant e+2\}\subset \Phi^{-1}(|K_{\sq}|)$ (in particular, if $z'$ is a vertex of the
boundary of $K$ in $\mbb T$ which is not one of the end points of the segment
$\widetilde L\cap K$, $d(z',\widetilde L)\geqslant e+2$). Let $K'_{\sq}$
be a compact subgraph of $\mbb T$ coinciding for $\sq=\alpha$ and $\sq=\beta$ such
that $\Phi_{\sq}^{-1}(|K'_{\sq}|)$ contains $U_{z,e+2}$ for every $z$ in
$\Phi^{-1}(K_{\sq})$.\\

Let $\Gamma=\Gal(\mbb T/\mbb G)$, let $H=\Stab(\widetilde L) (\simeq \mbf
Z)$ and let $\Gamma'$ be a subgroup of finite index of $\Gamma$ such that, for every $g\neq 1 \in \Gamma'$, $\Phi_{\sq}^{-1}(|K'|)\cap
g\cdot\Phi_{\sq}^{-1}(|K'|)=\emptyset$ and for every $g\in \Gamma\backslash H$,
$d_{\sq}(g\cdot\widetilde L,\widetilde L)>\diametre_{\sq}(|K'|)$($A=\{g\neq
1\in \Gamma|\Phi^{-1}(|K'|)\cap g\cdot\Phi^{-1}(|K'|)\neq \emptyset\}$ is finite by
compactness of $\Phi^{-1}(|K'|)$. So, as $\Gamma$ is residually finite, there
exists $\Gamma'_1$ of finite index in $\Gamma$ which does not intersect $A$. $B=\{g\in \Gamma/H-H |d_{\sq}(g\cdot \widetilde
L,\widetilde L)\leqslant \diametre_{\sq}(|K'|)\}$ is also finite, and as
 $\overline H\cap \Gamma=H$ in the profinite completion of $\Gamma$, there
exists $\Gamma'_2$ of finite index in $\Gamma$ containing $H$ such that
$\Gamma'_2\cap B\cdot H=\emptyset$. We may then chose $\Gamma'=\Gamma'_1\cap\Gamma'_2$). Let $H'=H\cap \Gamma'$. Let
$X'_{\sq}=\Omega_{\sq}/\Gamma'$ : its a finite topological covering of
$X_{\sq}$, and the isomorphism $\phi:\gtemp(X_\alpha)\simeq \gtemp(X_\beta)$
inducts an isomorphism
$\phi':\gtemp(X'_\alpha)\simeq\gtemp(X'_\beta)$.\\

Let $\mcal C_0$ be the current on $\mbb T$ with  $\mcal C_0(e)=+1$ if $e$ is an edge of
$\widetilde L$ (and $e$ has the good orientation) and $0$ otherwise (except if $e$
is an edge of $\widetilde L$ with the wrong orientation, in which case $\mcal C_0(e)=-1$) : this current is  invariant under $H$. Let $\mcal C_{\alpha}=\sum_{g\in
  \Gamma'/H'}g\cdot \mcal C_0$. It is a current on $\mbb T$, invariant under
$\Gamma'$ and which coincides with $\mcal C_0$ on $K'$.\\
Let $f_\alpha\in O(\Omega_\alpha)^*$ be the corresponding invertible function
with $f_\alpha(z_{0,\alpha})=1$, and let $X''_\alpha$ be a $(\mbf Z/p^e\mbf Z)$-torsor of $X'_\alpha$ corresponding to that current, that is to say such that its pullback $\widetilde X''_\alpha$ to
$\Omega_\alpha$ is isomorphic to $\Omega_\alpha\times_{\mbf G_m}\mbf G_m\to
\Omega_\alpha$ where the fiber product is taken, on the left side, along $f_\alpha$ and on the right side along $z\mapsto z^{p^e}$. Let
$X''_{\beta}=\phi'^*X''_\alpha$ (do not forget that $X''_\beta$ has no reason to correspond
to the current $\mcal C_\alpha$).\\
Let also $f_{0,\sq}\in O(\Omega_{\sq})^*$ be the invertible function corresponding to the current $\mcal C_0$ and let $\widetilde X_{0,\sq}$ be the corresponding $(\mbf
Z/p^e\mbf Z)$-torsor on $\Omega_{\sq}$. Recall that this
torsor is split over a point $z\in \mbb T(\Omega_{\sq})$ of $\Omega_{\sq}$
if and only if $d_{\sq}(z,\widetilde L)>e+\frac{1}{p-1}$.\\

According to corollary~(\ref{totdecrev}) applied to $U=\Phi^{-1}(|K'|)$ and to $V=\Phi^{-1}(|K|)$,
the torsor $\widetilde X''_\alpha-\widetilde X_{0,\alpha}$ on
$\Omega_\alpha$, which corresponds to the current $\mcal C_{\alpha}-\mcal C_0$ which
is zero over $K'$, is split over $\Phi^{-1}(|K|)$ since for every $z\in \Phi^{-1}(|K|)$, $U_{z,e+2}\in \Phi^{-1}(|K'|)$. Thus, for
$z\in K$, $\widetilde X''_\alpha$ is split if and only if $\widetilde X_{0,\alpha}$ is, if and only if $d_{\alpha}(z,\widetilde L)>e+\frac{1}{p-1}$.\\
In particular $\widetilde X''_\alpha$ is split over the vertices of the boundary of $K$ which are not the end points of $K\cap
\widetilde L$. Thus, according to the result of section~\ref{vertex}
applied to $X''_\alpha$ and $X''_\beta$ (as $\Omega_{\sq}\to
X'_{\sq}$ is a topological covering, $\widetilde X''_{\sq}\to
\Omega_{\sq}$ is split over a point if and only if
$X''_{\sq}\to X'_{\sq}$ over the image of that point), $\widetilde
X''_{0,\beta}$ is also split over the vertices of the boundary
of $K$ which are not the end points of $K\cap \widetilde L$.\\

Let $\mcal C_{\beta}$ be a current on $\mbb T(\Omega_\beta)$
corresponding to the $\mbf Z/p^e\mbf Z$-torsor $\widetilde X''_\beta$ (the current corresponding to
$\widetilde X''_\beta$  is well defined only modulo $p^e$). According to proposition~(\ref{nondecrev}), l
the restriction of $\mcal C_{\beta}$ to
the star of a vertex of $K$ which is not an end point of $K\cap\widetilde
L$ is zero modulo $p^e$. One deduces from it that, modulo $p^e$, the
restriction to $K$ of $\mathcal C_{\beta}$ must be congruent to the restriction
of $a\mcal C_0$. By adding to $\mcal C_{\beta}$ a current which is a multiple
of $p^e$, we may assume that $\mcal C_{\beta}-a\mcal C_0$ is zero on $K$ (because every current with boundary on $K$, that is to say that respect Kirchhoff's law in every vertex of the interior of $K$ but with no condition on the boundary of $K$, can be extended in a current on the whole $\mbb T$).\\
Thus, by applying corollary~(\ref{totdecrev}) to $U=\Phi^{-1}(|K|)$ and to
$V=\Phi^{-1}(K_0)$, one may deduce that $\widetilde X''_{\beta}-a\widetilde
X_{0,\beta}$ is split over $|K_0|$, so if $z$ is a vertex
of $K_0$, $\widetilde X''_{\beta}$ is split over
$z$ if and only if $\widetilde X_{0,\beta}$ is ($a$ is necessarily non zero modulo $p^e$ because $\widetilde X''_{\beta}$ cannot be split over $z'_0$), if and only if
$d_{\beta}(z,\widetilde L)>e+\frac{1}{p-1}$, according to lemma~(\ref{decompositionpuissance}).\\

Therefore, $d_{\beta}(z,\widetilde L)>e+\frac{1}{p-1}$ if and only if
$d_{\beta}(z,\widetilde L)>e+\frac{1}{p-1}$ if and only if
$d_{\alpha}(z,\widetilde L)>e+\frac{1}{p-1}$ for every vertex $z\in K_0$.\\
As one may choose $K_0$ as big as one wants, one may deduce that for every $z$ of $\mbb T$, $d_{\alpha}(z,\widetilde L)>e+\frac{1}{p-1}$
if and only if $d_{\beta}(z, \widetilde L)>e+\frac{1}{p-1}$, and this for every integer $e\geqslant 1$.\\
Thus $\max(1,\lceil d_{\beta}(z,\widetilde
L)-\frac{1}{p-1}\rceil)=\max(1,\lceil d_{\alpha}(z,\widetilde
L)-\frac{1}{p-1}\rceil)$.\\
By applying it to $(z_i)_{i\geqslant 0}$, one gets that for every
$i\geqslant 0$, \[\max(1,\lceil
i\lg_{\alpha}+r_{\alpha}-\frac{1}{p-1}\rceil)=\max(1,\lceil
i\lg_{\alpha}+r_{\alpha}-\frac{1}{p-1}\rceil).\]
One may deduce from this that for every loop $C$ of $\mbb G$ (and of every topological covering of $\mbb G$),
\[\lg_{\alpha}(C)=\lg_{\beta}(C).\]
One may conclude with the help of proposition~(\ref{combi}).\findem
\appendix
\section{A combinatorial result}
To end the proof of theorem~(\ref{mumford}), we will need to prove that one can recover from the length of all the loops of every finite covering of $\mbb G$ the length of every edge of $\mbb G$. We thus have to prove a general result on graphs such that the valency of every edge is at least 3.\\
To prove this, we will work by induction on the number of edges of the graph. However, by removing an edge of our graph, one does not in general get a graph such that the valency of every edge is at least three, and we may have to concatenate two edges or withdrawing an other edge to apply our induction assumption (but there will only be a few edges whose length cannot be recovered directly by the induction assumption). To exhibit enough loops that go through some edges we do not know yet the length thanks to the induction assumption, we will have to distinguish many cases depending on the number of connected components of the subgraph of all the edges whose length is already known.
\begin{prop}\label{combi} Let $\mbb G$ be finite graph such that the valency of every vertex is
  at least 3. Let $f:\{\text{edges of } \mbb G\}\to \mbf R$ be any function. Let
  us denote also by $f$ the induced function on the set of edge of a (topological) covering of
  $\mbb G$. Let us set, for $C$ a loop of a covering of
  $\mbb G$, $f(C)=\sum_{x\in\{\text{edges of } C\}}f(x)$.\\
If $f(C)=0$ for every loop $C$ of every covering of $\mbb G$, then
$f=0$.\end{prop}
\dem
Remark that if $\mbb G$ is a finite graph such that the valency of every vertex is at least 3
and if $\mbb H$ is a connected subgraph such that
the number of half-edges of $\mbb G\backslash \mbb H$ which ends in
$\mbb H$ is less than 3, then $\mbb H$ is not a tree
(if $\mbb H$ is a tree with at least one edge, $\mbb H$ has at
least two vertices of valency 1, and thus one already has 4 half-edges of $\mbb
G\backslash \mbb H$ which must end in one of those two vertices; if $\mbb
H$ is only a vertex, it is equally obvious).\\

Let us go back to our business. We will proceed by induction on the number of edges of $\mbb G$.\\
Thus let $(\mbb G,f)$ be a graph with $n\geqslant 1$ edges and a function $f$ on
the set of edges of $\mbb G$ which verify the hypotheses of the proposition, and
assume the proposition to be true if $\mbb G$ has less than $n$
edges.\\
We may assume $\mbb G$ to be connected (otherwise we may use the induction hypothesis to the various connected components to conclude).\\
Let $e$ be an edge of $\mbb G$, and start with assuming that:
\begin{enumerate}
\item the two end points of $e$ are different. Let $(m,n)$ be the arities of
  the end points of $e$ in $\mbb G\backslash\{e\}$. By our assumption about $\mbb
  G$, $m\geqslant   2$ and $n\geqslant 2$.
  \begin{enumerate} 
    \item If $m\geqslant 3$ and $n\geqslant 3$, $\mbb G'=\mbb
    G\backslash\{e\}$ is still a graph such that the valency of every vertex is at least $3$, thus one may apply our induction hypothesis to $\mbb G'$
    and to $f$ : $f(x)=0$, for every edge $\mbb G$ other than $e$.
    \begin{enumerate}
      \item If $\mbb G'$ is connected, then one may find a loop $C$ of $\mbb
        G$ that goes through $e$, and then $f(e)=f(C)=0$, which implies the result.
      \item If $\mbb G'$ has two connected components $A$ and $B$ (they cannot
        be trees according to the remark at the beginning of the proof), one may consider two coverings of order 2 			$A'$ and $B'$ of $A$ and $B$ respectively,
        that one patches in a covering $\mbb G'$ of $G$ of order 2.\\
        Then there exist a loop $C$ of $\mbb G'$ passing through the two edges over $e$.\\
        Then $2f(e)=f(C)=0$, which gives the wanted result.
    \end{enumerate}
    \item If $m\geqslant 3$ and $n= 2$ (or the other way), then let $a$
      and $b$ be the two edges starting at the second end point of $e$ (if $a$ and $b$
      are in fact the two half-edges of a single edge, the graph has then the same structure as the one of case 2.(c).i, where it will be treated; thus we will assume here that $a$ and $b$ are two different edges). Let $\mbb
      G'$ the graph obtained from $\mbb G$ by withdrawing $e$ and concatenating $a$ and $b$ in a single edge that 		  we will note $a+b$ (and we will define
      $f(a+b)=f(a)+f(b)$). Then $(\mbb G',f)$ satisfies the required conditions and so $f(x)=0$ for every edge of 		  $\mbb G$ other than $a,b$
      and $e$. Depending on the number of connected components of $\mbb G''=\mbb
      G\backslash\{a,b,e\}$, we will distinguish several cases:
     \begin{enumerate}
       \item 
           If $\mbb G''$ has only one connected component, then one may contract $\mbb
           G''$ to a point to get a graph $\mbb G_1$ with three edges
           (indeed, every loop $C_1$ of $\mbb G_1$ may be lifted
           to a loop $C$ of $\mbb G$ as $\mbb G''$ is connected, and
           $f(C_1)=f(C)$ since $f=0$ over $\mbb G''$), and thus $f(a)+f(b)=0$,
           $f(a)+f(e)=0$ and $f(b)+f(e)=0$, which gives the wanted result.
\begin{center}
\includegraphics{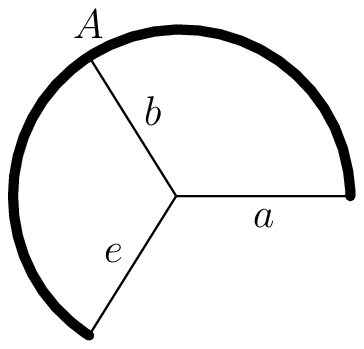}
\includegraphics{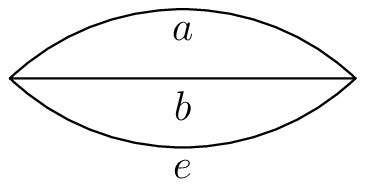}
\end{center}
       \item
           If $\mbb G''$ has two connected components $A$ and $B$ as
           in the picture (now $a$, $b$ and $e$ will play the same role
           and the proof will be the same if they are exchanged), then we start
           by considering a connected covering $\mbb G_1$ of order 2 of $\mbb G$
           such that its restrictions $A'$ and $B'$ to $A$ and to $B$ are connected (there exist such a covering
           as $A$ and $B$ cannot be trees according to the remark at
           the beginning of the proof), then one may
          contract $A'$ and $B'$ into a graph $G_2$. One gets
           $f(b)+f(e)=0$, $2f(a)+2f(b)=0$, $2f(a)+2f(e)=0$, which implies what we wanted.
 \begin{figure}[h]
 \begin{center}
 \includegraphics{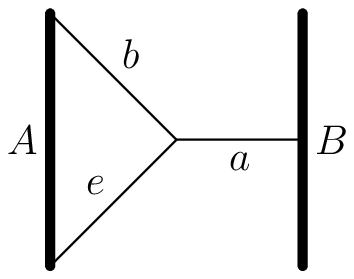}\hfill
 \includegraphics{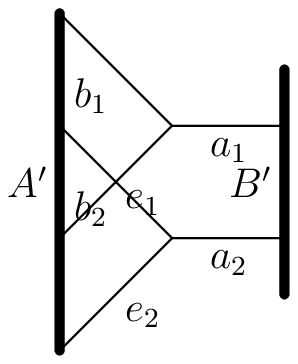}\hfill
 \includegraphics{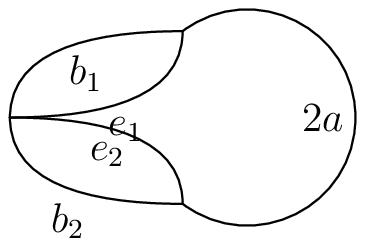}
 \end{center}
 \end{figure}
         \item
           If $\mbb G''$ has three connected components $A$, $B$ and
           $C$, then we start by considering a connected covering $\mbb G_1$
           of order 2 of $\mbb G$ such that its restrictions $A'$, $B'$ and $C'$ to
           $A$, $B$ and $C$ are connected, then one may contract $A'$, $B'$ and
           $C'$. One may deduce that $2f(a)+2f(b)=0,2f(b)+2f(e)=0,2f(e)+2f(a)=0$, which implies the result.
\begin{figure}[h]
\begin{center}
\includegraphics{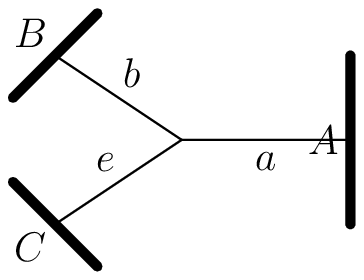}
\includegraphics{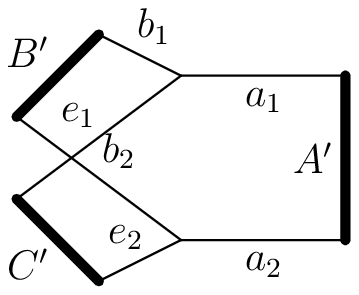}
\includegraphics{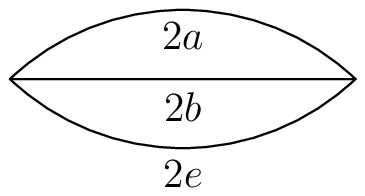}
\end{center}
\end{figure}
       \end{enumerate}
       \item If $m=n=2$, then let $a,b$ be the two (half-)edges ending at one end point of $e$ and  $c,d$ the
       (half-)edges ending at the other end point of $e$. Let $\mbb G'$ be the graph obtained
       from $\mbb G$ by withdrawing $e$ and concatenating $a$ and $b$ into
       $a+b$, and $c$ and $d$ into $c+d$. If one defines
       $f(a+b)=f(a)+f(b)$ and $f(c+d)=f(c)+f(d)$, $(G',f)$ satisfies
       assumptions of the proposition, and so, by the induction hypothesis,
       $f(x)=0$ for every edge $x$ of $\mbb G$ else than $a,b,c,d$ and
       $e$. Let $\mbb G''=\mbb G\backslash\{a,b,c,d,e\}$. Depending on the connected components of $\mbb G''$, we will distinguish several cases :
       \begin{enumerate}
         \item
           If $\mbb G$ has two connected components, $A$
           containing one end point of $a$ and of $b$, and $B$ containing one
           end point of $c$ and of $d$ (if $(a,b)$ or $(c,d)$ are the two
           half-edges of a single edge, the graph has the same structure as
           in case 2.(c).ii.B; if $(a,b)$ and $(c,d)$ both make single edges, then the structure is the one of the 			   degenerate case of 2.(c).(ii)), then start by considering a
           connected covering $\mbb G_1$ of order 2 of $\mbb G$ such that the
           restrictions $A'$ and $B'$ to $A$ and $B$ are connected, and then
           contract $A'$ and $B'$. One gets, for example, that
           $f(c)+f(d),f(a)+f(b),2(f(a)+f(e)+f(d)),2(f(a)+f(e)+f(c),2(f(b)+f(e)+f(d))=0$
           , which implies the result.
\begin{center}
\includegraphics{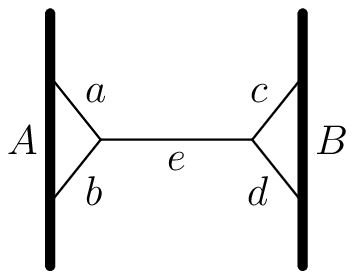}
\includegraphics{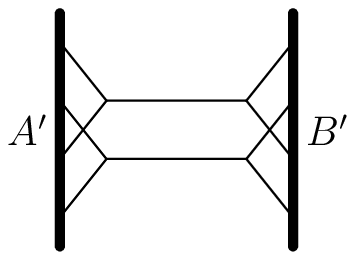}
\end{center}
           \item
             If $\mbb G$ has two connected components, $A$
             containing an end point of $a$ and $c$, and $B$ containing an
             end point of $b$ and $d$, then consider as usual a connected covering $\mbb G_1$ of order 2 of $\mbb
             G$ such that the restrictions $A'$ and $B'$ to $A$ and $B$ are
             connected, then contract $A'$ and $B'$. One gets for example that
             $f(a)+f(b)+f(c)+f(d)=2(f(a)+f(e)+f(d))=2(f(b)+f(e)+f(c))=2(f(a)+f(c))=f(b)+f(c)+f(e)=0$, which implies the result.
\begin{center}
\includegraphics{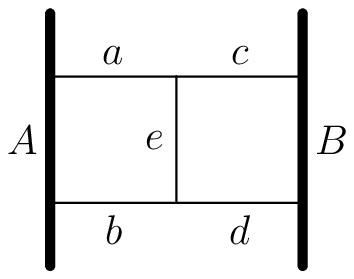}
\includegraphics{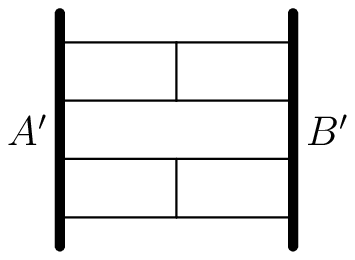}
\end{center}
             If $(a,c)$ (or symmetrically $(b,d)$) are the two
             half-edges of a single edge $a$, 
             consider a
             connected covering $\mbb G_1$ of order 2 of $\mbb G$ such that the
             restrictions $B'$ to $B$ and $(a\cup e)'$ and $a\cup e$ are connected. 
One gets that $f(a)+f(b)+f(c)= f(e)+f(b)+f(c)=
             2(f(a)+f(e)) = f(a)+f(e)+2f(b)=0$, which implies the result.
             \begin{center}
               \includegraphics{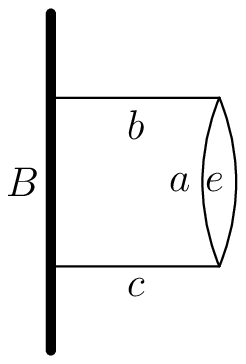}
               \includegraphics{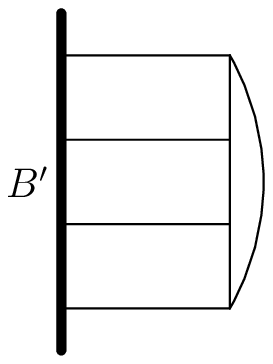}
             \end{center}
             If $(a,c)$ and $(b,c)$ are degenerated into two edges, then
             $\mbb G$ only has two vertices and three arrows joining those two vertices, and the result is obvious. 
             \item
             If $\mbb G$ has two connected components, $A$
             containing an end point of $a$,$b$ and of $c$, and $B$ containing an
             end point of $d$. We already know by induction hypothesis 
             applied to $\mbb G'$ that $f(c)+f(d)=0$. Now, by
             contracting $A$, one gets that
             $f(a)+f(b)=f(b)+f(c)+f(e)=f(a)+f(c)+f(e)=0$, which implies the result.
\begin{center}
\includegraphics{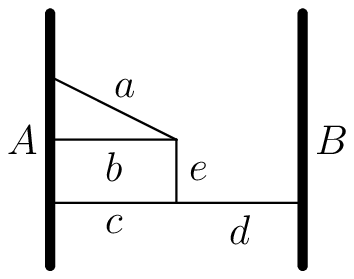}
\end{center}
             \item
             If $\mbb G$ has three connected components, $A$
             containing an end point of $a$, $B$ containing an end point of
             $b$, et $C$ containing an end point of $c$ and of $d$. Consider as usual a connected covering $\mbb G_1$ of order 2 of $\mbb
             G$ such that its restrictions $A'$, $B'$ and $C'$ to $A$, $B$ and $C$ are connected, then contract $A'$, $B'$ and $C'$. One gets for example that
             $f(c)+f(d) = 2(f(a)+f(b))=2(f(a)+f(e)+f(c)) =
             2(f(b)+f(e)+f(c)) = 2f(a)+2f(e)+f(c)+f(d)=0$, which implies the result.
\begin{center}
\includegraphics{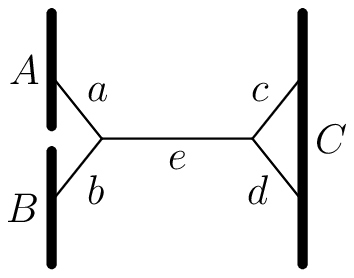}
\includegraphics{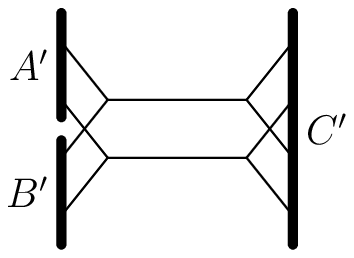}
\end{center}
             If $(c,d)$ is degenerated in a single edge $c$, the graph has the same structure as in case 2.(c).ii.A.
             \item
             If $\mbb G$ has three connected components, $A$
             containing an end point of $a$, $B$ containing an end point of
             $c$, and $C$ containing an end point of $b$ and $d$, then consider a connected covering $\mbb G_1$ of order 2 of $\mbb
             G$ such that its restrictions $A'$, $B'$ and $C'$ to $A$, $B$ and $C$ are
             connected then contract $A'$, $B'$ and $C'$. One gets for example that $f(b)+f(d)+f(e)=f(b)+f(d)+f(e)+2f(a)=2(f(a)+f(b))
             = 2(f(d)+f(c))=2(f(c)+f(e)+f(a))=0$.
             \begin{center}
                 \includegraphics{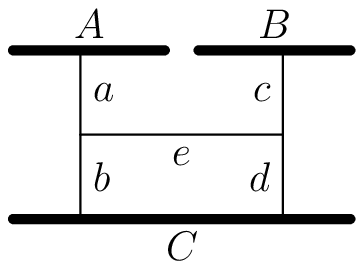}
                 \includegraphics{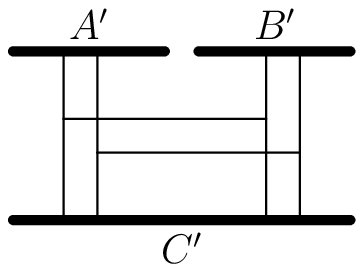}
               \end{center}
             If $(b,d)$ is degenerated into a single edge $b$, consider a
             covering $G_1$ of order 2 of $G$ such that its restrictions to $A$,
             $B$ and $b\cap e$ are connected, then contract $A'$ and $B'$. One gets that
             $2f(a)+f(e)+f(b) = 2(f(b)+f(e)) =2f(c)+f(b)+f(e) =
             2(f(a)+f(b)+f(c))=0$.
             \begin{center}
               \includegraphics{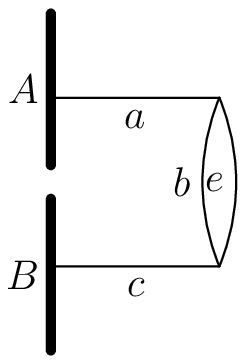}
               \includegraphics{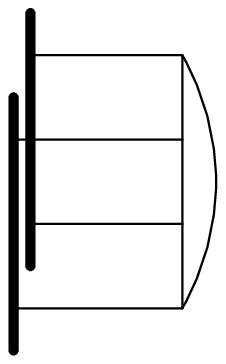}
             \end{center}
             \item
             If $\mbb G$ has four connected components, $A$
             containing an end point of $a$, $B$ containing an end point of $c$, $C$ containing an
             end point of $b$ and $D$ containing an end point of $d$, then consider
             a connected covering $\mbb G_1$ of order 2 of $\mbb
             G$ such that its restrictions $A'$, $B'$, $C'$ and $D'$ to $A$, $B$,
             $C$ and $D$ are
             connected then contract $A'$, $B'$, $C'$ and $D'$. One gets for example that
             $2(f(b)+f(d)+f(e))= f(b)+f(d)+f(e)+2f(a)=
             f(b)+f(d)+f(e)+2f(c)=2(f(b)+f(a))= 2(f(c)+f(d))=0$.
             \begin{center}
               \includegraphics{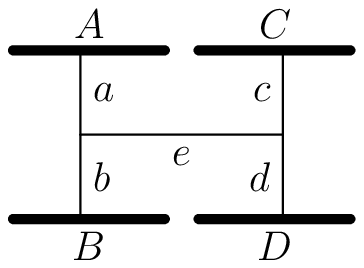}
               \includegraphics{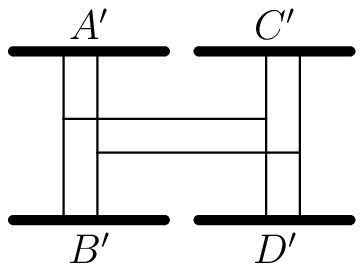}
             \end{center}
             \item
             If $\mbb G$ has a single connected component $A$,
             then contract it. One gets that $f(a)+f(b)=f(c)+f(d)=
             f(a)+f(c)+f(e)=f(b)+f(c)+f(e)=f(a)+f(d)+f(e)=0$, which implies the result.
             \begin{center}
               \includegraphics{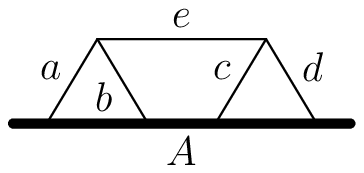}
             \end{center}
         \end{enumerate}
  \end{enumerate}
\item
the two end points of $e$ are a same vertex, and let $m$
be the valency of this vertex in $\mbb G\backslash\{e\}$. One has $m\geqslant 1$.
  \begin{enumerate}
    \item
      Assume $m\geqslant 3$. Then $\mbb G\backslash\{e\}$ satisfies the
      assumptions of the proposition and thus by induction hypothesis $f(x)=0$ for every edge of $\mbb G$
      else than $e$, and as $e$ is already a loop, $f(e)=0$ too.
    \item
      Assume $m=2$, and let $a$ and $b$ be the two edges ending at the end point of
      $e$ (if $a$ and $b$ are in fact only the two half-edges of a single edge,
      $\mbb G$ is only the wedge of two loops,
      and the result is obvious). Let $\mbb G'$ be the graph obtained from
      $\mbb G$ by withdrawing $e$ and by concatenating $a$ and $b$ in an
      edge $a+b$ and define $f(a+b)=f(a)+f(b)$. Then $\mbb G'$ satisfies the
      assumptions of the proposition and thus, by induction hypothesis, $f(x)=0$ 
      for every edge $x$ else than $e$, $a$ et $b$. Depending on the
      connected components of $\mbb G''=\mbb G\backslash\{a,b,e\}$, we will distinguish the following cases :
      \begin{enumerate}
        \item
          If $\mbb G''$ has a single connected component $A$, then consider a covering $\mbb
          G_1$ of order 2 of $\mbb G$ such that its restrictions to $A$ and $e$ are connected, then
          contract $A'$: one gets a graph $\mbb G_2$. One gets
          $2f(a)+f(e)= 2f(b)+f(e)=2f(e)=0$, which implies the result.
\begin{center}
\includegraphics{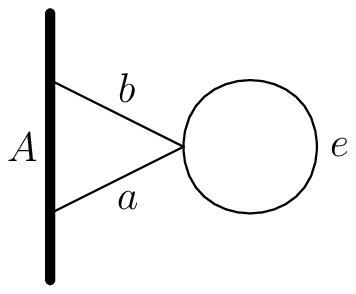}
\includegraphics{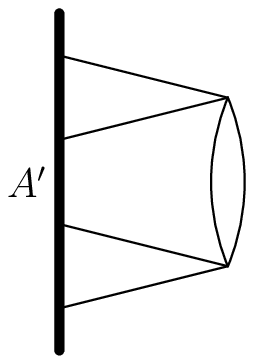}
\end{center}
        \item
          If $\mbb G''$ has two connected components, $A$ containing
          the end point of $a$, and $B$ containing the end point of $b$. Consider a covering $\mbb
          G_1$ of order 2 of $\mbb G$ such that its restrictions $A'$ and $B'$ to $A$ and $B$ are
          connected and such that its restriction to $e$ is connected too, then
          contract $A'$ and $B'$ to get a graph $\mbb G_2$. One gets that
          $2f(a)+f(e)= 2f(b)+f(e)=2f(e)=0$, which implies the result.
\begin{center}
\includegraphics{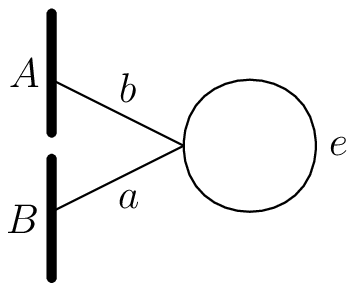}
\includegraphics{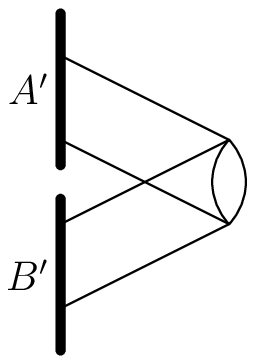}
\end{center}
      \end{enumerate}
   \item
     Assume $m=1$, and let $a$ be the edge ending at the end point of $e$. Let
     $n$ be the valency of the other end point of $a$. One must have
     $n\geqslant 2$.
     \begin{enumerate}
       \item
         Assume $n\geqslant 3$. Then $\mbb G'=\mbb G\backslash\{a,e\}$
         satisfies the assumptions of the proposition, so, by induction hypothesis, $f=0$ over $\mbb G'$. Consider now
         a covering $\mbb G_1$ of order 2 of $\mbb G$ such that its restrictions
         to $\mbb G'$ and to $e$ are connected, and contract
         the preimage of $\mbb G'$. One gets $2f(e)=2f(a)+f(e)=0$, which implies the result
\begin{center}
\includegraphics{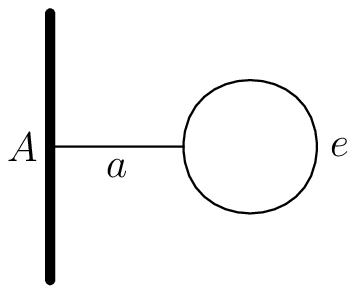}
\includegraphics{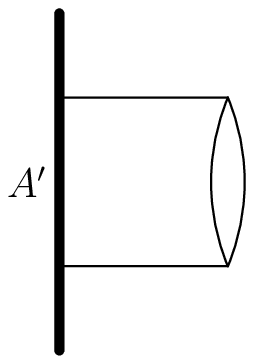}
\end{center}
       \item
         Assume $n\geqslant 2$, and let $c$ and $d$ be those two edges
         (if they only are the two half-edges of a single edge,
         $\mbb G$ is made of two loops joined by an edge; one shows the result for this particular graph by
         by considering the covering of order 2 whose restrictions to the loops are connected). Let $\mbb G''$ be the graph obtained from $\mbb G'$
         by concatenating $c$ and $d$ into $c+d$ and let us define
         $f(c+d)=f(c)+f(d)$. $\mbb G''$ satisfies the assumptions of the proposition,
         and so $f(x)=0$ for every edge x of $\mbb G$ else than $a,e,c$ and
         $d$ (and, in fact, $f(e)=0$ too). Let $\mbb G'''=\mbb G\backslash\{a,c,d,e\}$, and distinguish    the following cases depending on the number of connected components of $\mbb G'''$.
         \begin{enumerate}
           \item
             If $\mbb G'''$ has two connected components, $C$ containing
             the end point of $c$ and $D$ containing the end point of $d$,
             consider a covering $\mbb G_1$ of $\mbb G$ such that its restrictions 
             $C'$, $D'$ and $e'$ to $C$, $D$ and $e$ are connected, then contract
             $C'$, $D'$ and $e'$, one gets
             $2f(a)+2f(c)=2f(a)+2f(d)=2f(c)+2f(d)=0$, which implies the result.
             \begin{center}
               \includegraphics{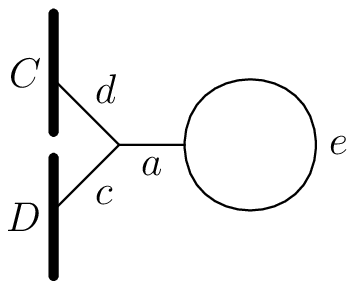}
               \includegraphics{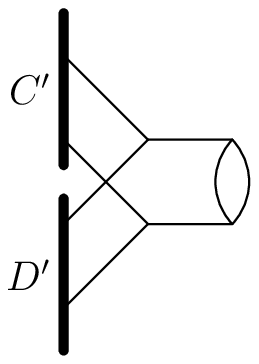}
             \end{center}
           \item
             If $\mbb G'''$ has a single connected component $A$, consider a
             covering $\mbb G_1$ of $\mbb G$ such that its restrictions $A'$
             and $e'$ to $A$ and $e$ are connected, then contract $A'$ and
             $e'$, one gets $2f(a)+2f(c)=2f(a)+2f(d)=2f(c)+2f(d)=0$, which implies the result.
             \begin{center}
               \includegraphics{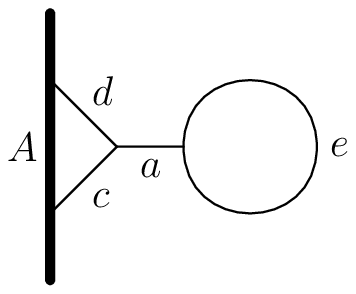}
               \includegraphics{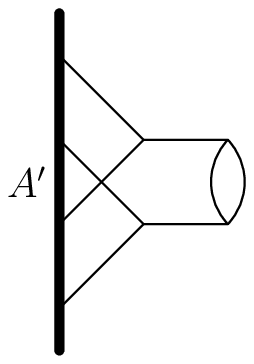}
               \end{center}
         \end{enumerate}
      \end{enumerate}
  \end{enumerate}
\end{enumerate} 
\findem

\bibliography{biblio2}
\bibliographystyle{amsplain}

\end{document}